\documentclass[pdflatex,sn-mathphys-num]{sn-jnl}

\usepackage{graphicx}%
\usepackage{multirow}%
\usepackage{amsmath,amssymb,amsfonts}%
\usepackage{amsthm}%
\usepackage{mathrsfs}%
\usepackage[title]{appendix}%
\usepackage{xcolor}%
\usepackage{textcomp}%
\usepackage{manyfoot}%
\usepackage{booktabs}%
\usepackage{algorithm}%
\usepackage{algorithmicx}%
\usepackage{algpseudocode}%
\usepackage{listings}%
\usepackage{comment}
\usepackage{enumerate}

\theoremstyle{thmstyleone}%
\newtheorem{theorem}{Theorem}[section]
\newtheorem{proposition}[theorem]{Proposition}
\newtheorem{lemma}[theorem]{Lemma}
\newtheorem{corollary}[theorem]{Corollary}

\theoremstyle{thmstyletwo}%
\theoremstyle{thmstylethree}%
\newtheorem{definition}{Definition}%
\newtheorem{remark}{Remark}%

\begin{document}

\title[Existence and nonexistence of global solutions to a semilinear parabolic equation with variable reaction on sub-Riemannian manifolds]{Existence and nonexistence of global solutions to a semilinear parabolic equation with variable reaction on sub-Riemannian manifolds}


\author*[1]{\fnm{Ricardo} \sur{Freire}}\email{ricardo.silva@uesb.edu.br}
\equalcont{These authors contributed equally to this work.}

\author[2]{\fnm{Ricardo} \sur{Castillo}}\email{rcastillo@ubiobio.cl}
\equalcont{These authors contributed equally to this work.}

\author[3]{\fnm{Omar}  \sur{Guzm\'an-Rea}}\email{oguzmanr@unmsm.edu.pe}
\equalcont{These authors contributed equally to this work.}

\author[4]{\fnm{Miguel} \sur{Loayza}}\email{miguel.loayza@ufpe.br}
\equalcont{These authors contributed equally to this work.}

\affil*[1]{\orgdiv{Departamento de Ciências Exatas e Tecnológicas}, \orgname{Universidade Estadual do Sudoeste da Bahia - UESB}, \orgaddress{\street{Estr. Bem Querer, Km-04}, \city{Vitória da Conquista}, \postcode{45083-900}, \state{Bahia}, \country{Brasil}}}

\affil[2]{\orgdiv{Facultad de Ciencias, Departamento de Matem\'atica}, \orgname{Universidad del B\'io-B\'io}, \orgaddress{\street{Avenida Collao 1202, Casilla 5-C}, \city{Concepción}, \postcode{4081112}, \state{B\'{\i}o-B\'{\i}o}, \country{Chile}}}

\affil[3]{\orgdiv{Facultad de Ciencias Matem\'aticas, Departamento de Matem\'atica}, \orgname{Universidad Nacional Mayor de San Marcos}, \orgaddress{\street{Ciudad Universitaria, Av. Venezuela s/n cuadra 34}, \city{Lima}, \country{Per\'u}}}

\affil[4]{\orgdiv{Departamento de Matem\'atica}, \orgname{Universidade Federal de Pernambuco}, \orgaddress{\street{Av. Prof. Moraes Rego, 1235-Cidade Universitária}, \city{Recife}, \postcode{50670-901}, \state{Pernambuco}, \country{Brasil}}}


\abstract{  We study the Cauchy problem $ u_t - \Delta_{\mathcal{M}} u = h(t)\,[f(u)]^{p(x)}$ on $\mathcal{M} \times (0,T)$, with $u(\cdot,0)=u_0 \in L^q(\mathcal{M})$, where $\mathcal{M}$ is a complete sub-Riemannian manifold endowed with a smooth volume measure and $\Delta_{\mathcal{M}}$ denotes the associated sub-Laplacian. The function $h \in L^1_{\mathrm{loc}}([0,\infty))$, $f$ is nonnegative and continuous, and $p$ is a bounded measurable function on $\mathcal{M}$. We establish sharp criteria for the global existence and nonexistence of nonnegative global solutions. In particular, our framework applies to connected unimodular Lie groups, yielding global existence, nonexistence, and Fujita-type results. Our approach allows for nonlinearities $f$ beyond the convex setting, thereby broadening the scope of the nonexistence results established in the Euclidean framework \cite[Theorem~1.5 (ii)]{CL}. 
Furthermore, our results generalize those of Ruzhansky and Yessirkegenov \cite[Theorems~1.1 and~1.3]{RY}, whose analysis relies on polynomial-type bounds of the form $K_1 u^p \le f(u) \le K_2 u^p$ with fixed constants $K_1,K_2>0$ and $p>1$.
}


	

\maketitle

\section{Introduction}

Let $T>0$. We consider the following semilinear parabolic problem
\begin{equation}\label{eq1}
 \left\{\begin{array}{rlll}
 u_t - \Delta_{\mathcal{M}} u  & = & h(t)[f(u)]^{p(\cdot)} & \mbox{ in } \mathcal{M} \times (0,T),\\
 \noalign{\medskip}
u(\cdot,0) &=&  u_{0} \geq 0   &     \mbox{ in } \mathcal{M}, \\
\noalign{\medskip}
u &\geq& 0  & \mbox{ in } \mathcal{M} \times (0,T),
      \end{array}\right.
\end{equation}
where $ 0 \leq u_0 \in L^q(\mathcal{M})$, $q>1$,  $h: (0,\infty) \rightarrow (0,\infty)$ belongs to $L^1_{\mathrm{loc}}([0,\infty))$, $f:[0,\infty)\to [0,\infty)$ is continuous, and $p$ is a bounded measurable function on $\mathcal{M}$ satisfying 
\begin{equation} \label{ExponentVar}1 < p^- := \operatorname*{ess\,inf}_{x\in\mathcal M} \,p(x)\leq p(x) \leq p^+ := \operatorname*{ess\,sup}_{x\in\mathcal M} \,p(x) < \infty.
\end{equation}
($p$ will be called a variable exponent). Here $\mathcal{M}$ is a sub-Riemannian manifold endowed with a smooth volume $\mu$, complete with respect to the associated sub-Riemannian distance, and  $\Delta_{\mathcal{M}}$ denotes the sub-Laplacian associated with $\mu$; further details will be provided in Section \ref{Prel}.

In the Euclidean setting $\mathcal{M}=(\mathbb{R}^N,+)$, problems involving variable exponents arise in the study of a wide range of nonlinear phenomena governed by elliptic and parabolic equations with non-standard diffusion. Such frameworks appear, for example, in mathematical models of non-Newtonian fluids with field-dependent constitutive laws, including electrorheological fluids whose mechanical properties depend on an applied electric field and thermoconvective flows with temperature-dependent viscosity \cite{Ruzika1, Antonsev1}, as well as in image restoration and denoising problems, where variable exponent diffusion models provide a flexible mechanism for adapting the smoothing process to local image features \cite{Aboulaich, Bollt}. 
Moreover, variable exponent spaces are central to the study of elliptic and fluid-mechanics problems with non-standard growth, such as the Poisson and Stokes equations \cite[Chapter~14]{DIENING}. 
They also provide an appropriate functional framework for the study of existence and regularity of parabolic equations with variable growth; see, for example, \cite{Arora} and references therein.

\subsection*{Euclidean Framework}
The study of sharp criteria for global existence and nonexistence of global solutions for equation~(\ref{eq1}), in the case where $h(t) \equiv 1$, $f(u)=u$, and $p(x)\equiv p$ is constant, with the abelian structure $\mathcal{M}=(\mathbb{R}^N,+)$, originates from the classical work of Fujita~\cite{Fujita}, who introduced the critical exponent
$$p_c= 1 + \frac{2}{N},$$
which separates the values of $p$ for which global solutions exist from those for which no global solutions exist.
More precisely, if $1 < p \leq p_c$, then all nontrivial nonnegative solutions blow up in finite time, whereas for $p > p_c$ there exist global solutions for sufficiently small initial data. 
Since then, these seminal ideas have been developed in various directions by numerous researchers; see, for example, \cite{quittner,AronsonW,Weissler,LeeNi,meier3,Ferreira,Castillo,Carhuas} and the references therein.

Building upon earlier contributions, Ferreira, de Pablo, Pérez-Llanos, and Rossi \cite{Ferreira} analyzed problem \eqref{eq1} in the case where $h(t) \equiv 1$, $f(u)=u$, and $p(x)$ is bounded and continuous satisfying \eqref{ExponentVar}, in $\mathbb{R}^N$ and in bounded smooth domains. They proved that, in $\mathbb{R}^N$, all nontrivial nonnegative solutions blow up in finite time when $1<p^{-}<p^{+}<p_c$, whereas global solutions exist when $p^{-}>p_c$. In the case $p^{-}<p_c<p^{+}$, they showed that there exist functions $p(x)$ for which problem \eqref{eq1} admits global nontrivial solutions and others for which all solutions blow up in finite time. In particular, their work extends Fujita’s classical results \cite{Fujita} to the framework of variable exponents $p(x)$. 
Castillo and Loayza \cite{Castillo} subsequently extended these results by employing differential inequalities and semigroup techniques, 
which allowed them to address the global existence and nonexistence of global solutions to problem \eqref{eq1} on arbitrary domains, requiring only information on the asymptotic behavior of the associated Dirichlet semigroup. For further developments in this direction, see, for instance, \cite{meier2,meier3,Castillo1,Castillo2,Castillo3,Castillo4,Carhuas}. 
A gap remained in \cite{Castillo} 
regarding the nonexistence of global solutions and the possible existence of unbounded global solutions. This issue was later resolved by the same authors in \cite{CL} through the analysis of a more general nonlinearity $f$ satisfying, in particular, Jensen's inequality \eqref{Jensen11}. 

Many nonexistence results for global solutions (or local nonexistence results) available in the literature rely on Jensen's inequality
\begin{equation} \label{Jensen11} f(S(t)u_0) \leq S(t)f(u_0),
\end{equation}
which is satisfied, for example, when $f$ is convex and $f(0)=0$. In this context, $\{S(t)\}_{t\geq0}$ denotes the Dirichlet heat semigroup. For comprehensive discussions and further developments, see \cite[Theorem 1.5]{Miguel1}, \cite[Theorem 2 - (ii)]{CastilloDeg}, \cite[Theorem 1.5 - (ii)]{CL}, \cite[Theorem 6 - (ii)]{Carhuas1}, \cite[Theorem 1.2]{Kumar1}, \cite[Lemma 4.4]{Lenz2023}, \cite[Proof of Theorem 2(ii)]{Castillo2026}, \cite[Remark 2.11]{Meglioli2026}, and the references therein.  See also \cite[Lemma 4.1]{Fujishima1}, where Fujishima and Ioku derived a necessary condition for the existence of solutions and used it to obtain nonexistence results for certain classes of initial data.

The preceding discussion naturally raises the question of whether nonexistence results for global solutions can be established
without assuming Jensen's condition \eqref{Jensen11}. One of the objectives of the present work is to address this issue. In Theorem \ref{thm1} (ii), we establish
a nonexistence result for global solutions
without requiring \eqref{Jensen11}. Our approach is developed in a flexible framework that encompasses diverse geometric structures.

\subsection*{The Non-Euclidean Setting}

The study of global existence, nonexistence of global solutions, and Fujita-type results in non-Euclidean geometries and subelliptic operators began to emerge in the 1990s. During that decade, in the case where $h(t)\equiv 1$, $f(u)=u$, and $p(x)\equiv p>1$, Pascucci \cite{Pascucci1} studied problem \eqref{eq1} for the sub-Laplacian on nilpotent stratified Lie groups over $\mathbb{R}^N$, and later in \cite{Pascucci1999} established Fujita-type results for a class of hypoelliptic operators of degenerate parabolic type, including a linearized prototype of the Kolmogorov operator. Around that time, Zhang \cite{Zhang3} investigated the problem for general subelliptic operators, including the sub-Laplacian on the Heisenberg group, while in \cite{Zhang1, Zhang2} he considered the Laplace--Beltrami operator on noncompact Riemannian manifolds $\mathcal{M}$ with $\dim(\mathcal{M})\geq 3$. In \cite{Pascucci1,Pascucci1999,Zhang3}, the Fujita exponent
$$
p_c = 1 + \frac{2}{Q},
$$
depends on the homogeneous dimension $Q$ of the stratified Lie group $(\mathbb{R}^{N}, \circ)=\mathcal{M}$ with group law $\circ$, which is generally strictly larger than the topological dimension; for example, in the Heisenberg group $\mathcal{M}=\mathbb{H}^{N}=\mathbb{R}^{2N+1}$ one has $Q=2N+2$. Moreover, the homogeneous dimension $Q$ is closely related to the polynomial volume growth of the group. On the other hand, in \cite{Zhang2}, Zhang assumed nonnegative Ricci curvature and that the volume of geodesic balls in $\mathcal{M}$ grows like $R^{\alpha}$ as $R\to\infty$, where $R$ denotes the radius of the balls and $\alpha>2$. In this framework, the Fujita critical exponent is given by
$$
p_c=\frac{\alpha}{\alpha-2}.
$$
Subsequently, the same author extended these results to weakly coupled parabolic systems \cite{Zhang4} and to problems involving more general spatially weighted source terms and higher-order operators.

The preceding results show that the Fujita critical exponent in the non-Euclidean frame is strongly connected with the asymptotic polynomial volume growth of geodesic balls in the underlying space. In particular, for stratified Lie groups, this behavior is encoded by the homogeneous dimension $Q$, whereas on noncompact Riemannian manifolds, it is determined by the large-scale geometry. 

A substantially different behavior arises in the hyperbolic setting $\mathbf{H}^{N}$ due to the exponential volume growth of geodesic balls, which falls outside the assumptions considered in \cite{Zhang2}. In this direction, Bandle, Pozio, and Tesei~\cite{CBandle1} studied problem~\eqref{eq1} (with $h(t)=e^{\alpha t}$ for $\alpha>0$, $f(u)=u$, and $p(x)=p>1$) on the hyperbolic space $\mathbf{H}^{N}=\mathcal{M}$ and proved that
$$p^{*}_{\mathbf{H}^{N}} = 1 + \frac{\alpha}{\lambda_0},
   \qquad \text{where }\lambda_0= \frac{(N-1)^2}{4},$$
is critical in the Fujita sense when $p^{*}_{\mathbf{H}^{N}}\neq p$. In this case, the Fujita exponent for the hyperbolic space coincides with that for bounded domains in the Euclidean setting. By contrast, in the critical case $p^{*}_{\mathbf{H}^{N}}=p$, global solutions exist, whereas in bounded Euclidean domains every positive solution blows up in finite time (see \cite[Appendix A and Remark 2]{Wang2}).

The above results have been further developed in several directions by many authors; see, for example, \cite{Suragan1,HCastillo,Freire1,Georgiev1,Grillo1,Punzo1,Wang1,Monticelli1} and the references therein. 
More recently, Ruzhansky and Yessirkegenov~\cite{RY} established results on the global existence and nonexistence of solutions to problem \eqref{eq1} with $h(t)\equiv 1$ and $p(x)=1$ on sub-Riemannian manifolds, a setting which encompasses both Riemannian manifolds and stratified Lie groups in a unified framework. In particular, in the case of connected unimodular Lie groups with polynomial volume growth, their analysis yields the corresponding Fujita critical exponent. Specifically, they proved the following results:

\begin{theorem}[\cite{RY}] \label{Rush1} \begin{enumerate}
\item[(i)] Consider a sub-Riemannian manifold $\mathcal{M}$. Suppose $f:[0,\infty)\to[0,\infty)$ is a continuous increasing function satisfying
$$f(u)\leq K_1\, u^q, \quad \mbox{for all } u\in[0,\infty),$$
for some constant $K_1>0$ and exponent $1<q<\infty$. Let $u_0\in L^r(\mathcal{M})$ be nonnegative, with $1\leq r<\infty$, and assume that
\begin{equation*}
\int_0^\infty \bigl\| S_{\mathcal{M}}(s)u_0 \bigr\|_{L^\infty(\mathcal{M})}^{\,q-1}\, ds
< \frac{1}{K_1(q-1)} .
\end{equation*}
($\{ S_{\mathcal{M}}(t) \}_{t\geq 0}$ denotes the heat semigroup associated with $\mathcal{M}$). Then there exists a nonnegative continuous curve $ 
u:[0,\infty)\to L^r(\mathcal{M})
$ 
which is a global solution to \eqref{eq1} with initial condition $u_0$. 
\item[(ii)] Let $\mathcal{M}=\mathcal{G}$ be a connected unimodular Lie group whose volume growth is polynomial of order $D$, and let $1<q<\infty$. Denote
$q_F := 1+\frac{2}{D}.$
\begin{enumerate}
\item[(a)]
Assume that $q\leq q_F$ and let $f:[0,\infty)\to[0,\infty)$ be a locally integrable function satisfying
$$
f(u)\ge K_2 u^q, \quad \mbox{for all } u\in[0,\infty),
$$
for some constant $K_2>0$. Then \eqref{eq1} admits no nontrivial nonnegative global solution.
\item[(b)]
If $q>q_F$ and $f:[0,\infty)\to[0,\infty)$ is a continuous increasing function such that
$$
f(u)\le K_1 u^q, \quad \mbox{for all } u\in[0,\infty),
$$
for some constant $K_1>0$, then for every $1\le r<\infty$ the associated Cauchy problem \eqref{eq1} possesses a nonnegative global solution corresponding to some nonnegative initial datum
$u_0\in L^r(\mathcal{G})$.
\end{enumerate}
\end{enumerate}
\end{theorem}
In addition, the authors proved global existence results when $\mathcal{G}$ is a connected unimodular Lie group with exponential volume growth.

Motivated by the results discussed above, a further objective of the present work is to investigate problem \eqref{eq1} in the setting of sub-Riemannian manifolds while relaxing the global growth assumptions
$$ K_2 u^q \leq f(u) \leq K_1 u^q, \quad \mbox{for all } u\in[0,\infty) $$
required in Theorem \ref{Rush1} (see Theorem~\ref{thm1}).

\subsection*{Organization of the paper}

 The remainder of the paper is organized as follows. Section \ref{Prel} contains the necessary preliminaries, while Section \ref{mainnnn} presents the main results. Their proofs are given in Section \ref{proooofs}. In Section~\ref{connlieg}, we apply our results to connected unimodular Lie groups with either polynomial or exponential volume growth, thereby covering a broad class of nonlinearities $f$ (see Theorems~\ref{apli-lie}, \ref{lieee22222}, \ref{apli-lie-case-exp}, \ref{lieee222223} and \ref{apli-lie-case-conv}). 
For example, our framework includes nonconvex functions that do not satisfy the polynomial lower bound condition (see Theorem \ref{log-exam}). In particular, for $f(u)=u$ and $h(t)=t^q$ with $q>-1$, on connected unimodular Lie group whose volume growth is polynomial of order $D$,
we prove that global solutions exist whenever 
$$p^{-} > 1 + \frac{2(q+1)}{D},$$
whereas solutions cannot exist globally if
$$p^{+} < 1 + \frac{2(q+1)}{D}.$$
Finally, Section~\ref{SectionHeisen} presents further results, including an analysis of the intermediate regime
$$p^{-} \leq 1 + \frac{2(q+1)}{D} \leq p^{+}$$
for the case $f(u)=u$ and $h(t)=t^q$. On homogeneous Carnot groups, we identify a class of initial data for which no global solution exists, independently of the variable exponent $p(x)$ (see Corollary \ref{co-caso-crit1}). On the other hand, we exhibit another class of initial data for which global solutions exist whenever
$$p^{-} > 1 + \frac{2(q+1)}{D}$$
(see Corollary~\ref{co-caso-crit2}), and we identify a second Fujita critical exponent in the sense of Lee and Ni \cite{LeeNi} (see Corollary~\ref{co-caso-crit3}). Throughout the paper, $C, C_1, C_2, \dots$ denote positive constants that may vary from line to line without affecting the validity of the results.

\section{Preliminaries}\label{Prel}

We begin this section by recalling some fundamental aspects of sub-Riemannian geometry, as can be found in detail in \cite{ABB, barl, ovi, stri1}.

A sub-Riemannian manifold is a triple $(\mathcal{M}, \mathcal{H}, g)$, where: $\mathcal{M}$ is a $N$-dimensional connected orientable smooth manifold; $\mathcal{H}$ is a smooth subbundle of the tangent bundle $T\mathcal{M}$ (called horizontal distribution) of constant rank $k < N$ satisfying the H\"ormander condition (also known as the bracket-generating condition, see \cite[Definition 2.7.1]{ovi} or \cite[Section 2]{stri1}). This means that the Lie algebra generated by smooth sections of $\mathcal{H}$ spans the entire tangent space $T_q \mathcal{M}$ at each $q \in \mathcal{M}$. The sub-Riemannian metric $g$ is given by $g_q:\mathcal{H}_q \times \mathcal{H}_q \to \mathbb{R}$ for each $q \in \mathcal{M}$, where $g_q$ is a positive-definite quadratic form on $\mathcal{H}_q$, the fiber of $\mathcal{H}$ at $q$, and $g_q$ depends smoothly on $q$. 

A Lipschitz continuous curve $\gamma : [0, T] \to \mathcal{M}$ is said to be horizontal if 
$$
\dot{\gamma}(t) \in \mathcal{H}_{\gamma(t)} \quad \text{for almost every } t \in [0, T].
$$
Given a horizontal curve $\gamma : [0, T] \to \mathcal{M}$, the length of $\gamma$ is defined as
\begin{equation*}
    l(\gamma) = \int_0^T \sqrt{g( \dot{\gamma}(t),\dot{\gamma}(t) )} \, dt.
\end{equation*}
The distance induced by the sub-Riemannian structure on $\mathcal{M}$ is the function
\begin{equation*}
    d_C(q_0, q_1) = \inf\{l(\gamma) \mid \gamma(0) = q_0, \ \gamma(T) = q_1, \ \gamma \text{ horizontal curve}\}.
\end{equation*}

The hypothesis of connectedness of $\mathcal{M}$ and the H\"ormander condition guarantees the finiteness of $d_C$ with respect to the topology of $\mathcal{M}$, as stated by Chow's Theorem (see \cite[Theorem 3.5.2]{ovi}). The function $d_C$ is called the Carnot-Carath\'eodory distance of the sub-Riemannian manifold $\mathcal{M}.$

The sections in the horizontal bundle $(\mathcal{H} , \mathcal{M})$ are called horizontal vector fields. A system of $k$ linearly independent horizontal vector fields  $\{X_1, \ldots, X_k\}$ that span $\mathcal{H}$ and are orthonormal with respect to $g$, i.e.,
\begin{equation*}
    \mathcal{H}_q = \operatorname{span}\{X_1(q), \ldots, X_k(q)\}, \quad g_q(X_i(q), X_j(q)) = \delta_{ij},
\end{equation*}
 is called a local orthonormal frame. The existence of such a system can be found in \cite[Theorem 2.3.4]{ovi}. The horizontal gradient on $\mathcal{M}$ with respect to the local orthonormal frame  $\{X_1, \ldots, X_k\}$ is defined as 
$$
\nabla \varphi = \sum_{i=1}^k X_i(\varphi) X_i,
$$
where $X_i(\varphi)$ is the Lie derivative of $\varphi$ in the direction of $X_i$. 

The divergence $\operatorname{div}_\mu X$ of a vector field $X$ with respect to a smooth volume measure $\mu$, is defined by the identity $L_X \mu = (\operatorname{div}_\mu X) \mu,$ where $L_X$ is the Lie derivative with respect to $X$. The sub-Laplacian associated with the sub-Riemannian manifold is defined as $\Delta_{\mathcal{M}} \varphi = \operatorname{div}_\mu (\nabla \varphi),$
and it can be written in a local orthonormal frame $X_1, \ldots, X_k$ as
\begin{equation*}
\Delta_{\mathcal{M}} = \sum_{i=1}^k X_i^2 + (\operatorname{div}_\mu X_i) X_i.
\end{equation*}
Therefore, the sub-Laplacian $\Delta_{\mathcal{M}}$ is the natural generalization of the Laplace-Beltrami operator defined on a Riemannian manifold. 

In general, the sub-Laplacian $\Delta_{\mathcal{M}}$ can always be expressed as the sum of squares of the elements of the orthonormal frame plus a first-order term that depends on the choice of the volume $\mu$. Moreover, up to an adjustment in the first-order term, $\Delta_{\mathcal{M}}$ is formally symmetric with respect to $L^2(\mathcal{M}, \mu)$ (see \cite[Section 12]{stri1}). The results of Strichartz \cite{stri2} and \cite{stri1}, guarantee that if the manifold $\mathcal{M}$ is complete in the sub-Riemannian metric associated with $\Delta_{\mathcal{M}}$ and provided smooth volume $\mu$, then $\Delta_{\mathcal{M}}$ is essentially self-adjoint on $C^{\infty}_c(\mathcal{M})$. Therefore, if we consider $\{S_{\mathcal{M}}(t)\}_{t \geq 0}$ to be the heat semigroup on $L^2(\mathcal{M}, \mu)$ defined by the spectral theorem, then there exists a symmetric, positive  $C^{\infty}$ kernel $H_t(x, y)$ on $(0, \infty) \times \mathcal{M} \times \mathcal{M}$ such that 
\begin{equation}\label{energ-kernel}
    \int_{\mathcal{M}} H_t(x, y) d \mu (y) \leq 1
\end{equation}
for all $x\in \mathcal{M}$ and $t>0$, and
\begin{equation}\label{conv-semigroup}
    [S_{\mathcal{M}}(t)u](x) = \int_{\mathcal{M}} H_t(x, y) u(y) \, d\mu(y)
\end{equation} 
for all $u \in L^2(\mathcal{M}, \mu)$. Moreover, \eqref{conv-semigroup} satisfies the sub-Riemannian heat equation, that is,
$$\partial_t(S_{\mathcal{M}}(t)u) - \Delta_{\mathcal{M}}(S(t)_{\mathcal{M}}u) = 0$$ 
for every $u \in L^2(\mathcal{M}, \mu)$. 

Defining $S_{\mathcal{M}}(t) u$ by \eqref{conv-semigroup} for all $u \in L^q(\mathcal{M}, \mu)$, $1 \leq  q \leq  \infty$, it follows that $S_{\mathcal{M}}(t) u$  remains a solution to the sub-Riemannian heat equation. Furthermore, the semigroup $\{S_{\mathcal{M}}(t)\}_{t \geq 0}$ is strongly continuous on
$L^q(\mathcal{M}, \mu)$, for $1 \leq  q <  \infty$ and contractive for all $L^q(\mathcal{M}, \mu)$, $1 \leq  q \leq  \infty$ (for further details, see \cite[Section 3]{stri2} and \cite[Section 12]{stri1}).

  Henceforth, $\mathcal M$ denotes a complete sub-Riemannian manifold endowed with a smooth $\sigma$-finite volume measure $\mu$. The heat semigroup $\{S_{\mathcal M}(t)\}_{t \ge 0}$ on $L^q(\mathcal M,\mu)$, $1 \le q \le \infty$, associated with the sub-Laplacian $\Delta_{\mathcal M}$ is strongly continuous on
$L^q(\mathcal{M}, \mu)$, for $1 \leq  q <  \infty$ and contractive for all $L^q(\mathcal{M}, \mu)$, $1 \leq  q \leq  \infty$, and is understood in terms of convolution with the heat kernel, as in \eqref{conv-semigroup}. We now define what we mean by a solution to problem \eqref{eq1}.

\begin{definition}\label{Def}
Let $1 \leq r \leq \infty$ and let $u_0\in L^{r}\left ( \mathcal{M} , \mu \right ) $ with $u_0\geq 0$. Let $h, \,f\in C([0,\infty))$ and assume furthermore that $p : \mathcal{M} \to [1, \infty)$ is a bounded measurable function satisfying $1 < p^{-} \leq p(x) \leq p^{+}  < + \infty$ for almost every (a.e.) $x \in \mathcal{M}$. A measurable nonnegative function $u$ defined a.e. in $\mathcal{M} \times \left ( 0,T \right )$, is said to be a solution to problem (\ref{eq1}) if $u\in L^{\infty }\left ( \left ( 0,T \right ),   L^{r}\left ( \mathcal{M} , \mu \right )  \right )$ and 
\begin{equation}\label{eq2}
   u\left ( t \right )=S_{\mathcal{M}}\left ( t \right ) u_0+\int_{0}^{t}S_{\mathcal{M}}\left ( t-\sigma \right ) h(\sigma) [f\left ( u\left ( \sigma \right ) \right )]^{p(\cdot)}d\sigma
\end{equation}
a.e. in $\mathcal{M} \times \left ( 0,T \right )$, where the action of the operator $S_{\mathcal{M}}(t)$ is given by the integral representation \eqref{conv-semigroup}. When $T=\infty$, the solution is said to be global.   
\end{definition}
\begin{remark} The $\sigma$-finiteness of the measure $\mu$ plays an important role. In particular, it is required to ensure the completeness of the variable exponent Lebesgue space $L^{p(x)}(\mathcal{M}, d\mu)$ and to allow the application of Tonelli’s theorem, as noted in Remark~\ref{remarks-princteo}(c). In the case of an equiregular sub-Riemannian manifold, the measure $\mu$ may be chosen as Popp’s volume, see
\cite{bar11}.
\end{remark}

\section{Main Results}\label{mainnnn}

Before stating our main result, we establish some technical results.
\begin{lemma} Let $\eta$ be a positive measure in $\mathcal{M}$ such that $\textstyle\int_{\mathcal{M}} d\eta \leq 1$ and let $f\in L^{p^+}(\mathcal{M}, d\eta)$ with $1\leq p^{-} \leq p(x)\leq p^{+} < + \infty$ for almost every $x \in \mathcal{M}$. Then
\begin{equation*}
\int_{\mathcal{M}} \vert f(x)\vert ^{p(x)} d\eta(x) \geq 2^{-p^+} \min\left\{ \left( \int_{\mathcal{M}} \vert f(x)\vert d\eta(x) \right)^{p^-}, \left(\int_{\mathcal{M}} \vert f(x)\vert d\eta(x) \right)^{p^{+}} \right\}.    
\end{equation*}
\label{Ferreira-l2}
\end{lemma}

In the Euclidean setting, Lemma \ref{Ferreira-l2} was established in \cite[Lemma 3.1]{Ferreira}. Before proceeding to the proof of Lemma \ref{Ferreira-l2}, we provide a brief overview of Lebesgue spaces with variable exponents on manifolds. Regarding the manifold $\mathcal{M}$ as a complete $\sigma$-finite measure
space, we consider a bounded measurable function
$p : \mathcal{M} \to [1,\infty)$ and define
$$
p^{-} := \operatorname{ess\,inf}_{x \in \mathcal{M}} p(x), \quad p^{+} := \operatorname{ess\,sup}_{x \in \mathcal{M}} p(x).
$$

The variable exponent Lebesgue space $L^{p(x)}(\mathcal{M}, 
d\eta)$,  which is a particular case of a Musielak--Orlicz space (cf.~\cite[Section 3.2]{DIENING}), is defined in our setting as the set of all measurable functions $g: \mathcal{M} \to \mathbb{R}$ such that
$ \textstyle\int_{\mathcal{M}} \vert g(x)\vert ^{p(x)} d\eta(x) < \infty$. This space becomes a Banach space when endowed with the Luxemburg norm, defined by
$$ \Vert g\Vert _{L^{p(x)}(\mathcal{M}, d\eta)} = \inf \left\{ \tau > 0 : \int_{\mathcal{M}} \left\vert  \frac{g(x)}{\tau} \right\vert ^{p(x)} d\eta(x) \leq 1 \right\}. $$
(see \cite[Theorem 3.2.7]{DIENING}).
In particular, we have the Hölder inequality \cite[Lemma 3.2.20]{DIENING}
\begin{equation}\label{holder-variavel}
    \Vert uv\Vert _{ L^{r(x)}(\mathcal{M}, d\eta)} \leq 2 \Vert u\Vert _{L^{p(x)}(\mathcal{M}, d\eta)} \Vert v\Vert _{L^{q(x)}(\mathcal{M}, d\eta)}
\end{equation}
provided that $1/r(x) = 1/p(x) + 1/q(x)$ holds for almost every $x \in \mathcal{M}$. For further aspects regarding these variable exponent spaces, such as reflexivity, variable exponent Sobolev spaces, and embeddings, we refer the reader to \cite{KOV}, \cite{FAN}, and \cite{DIENING}.

\begin{proof}[Proof of Lemma \ref{Ferreira-l2}]
Following the arguments in \cite[Lemma 3.1]{Ferreira}, we observe that $f\in L^{p^+}(\mathcal{M}, d\eta)$ ensures that $f \in L^{p(x)}(\mathcal{M}, d\eta)$. For any nonzero function $f$, it can be shown that
$$ \int_{\mathcal{M}} \left\vert  \frac{f(x)}{\|f\|_{L^{p(x)}(\mathcal{M}, d\eta)}} \right\vert ^{p(x)} d\eta(x) =1, $$
from which it follows that
$$ \int_{\mathcal{M}} \vert f (x)\vert ^{p(x)} d\eta(x) \geq \begin{cases} \|f\|_{L^{p(x)}(\mathcal{M}, d\eta)}^{p^{-}} & \text{if } \|f\|_{L^{p(x)}(\mathcal{M}, d\eta)} \geq 1, \\
\noalign{\medskip}
\|f\|_{L^{p(x)}(\mathcal{M}, d\eta)}^{p^+} & \text{if } \|f\|_{L^{p(x)}(\mathcal{M}, d\eta)} \leq  1. \end{cases} $$
On the other hand, by Hölder's inequality \eqref{holder-variavel}, it follows that
$$ \int_{\mathcal{M}} \vert f(x)\vert  d\eta(x) = \int_{\mathcal{M}} \vert f(x)\vert  \cdot 1 \,d\eta(x)  \leq 2 \|f\|_{L^{p(x)}(\mathcal{M}, d\eta)}, $$
where we have used the fact that $\textstyle\int_{\mathcal{M}} d\eta \leq 1$. This completes the proof.
\end{proof}

We shall establish an important version of Jensen's inequality, which will be needed
in the proof of Theorem~\ref{thm1} (ii).

\begin{lemma}\label{jensen}  If $\,0\leq u \in L^{\infty}(\mathcal{M}, \mu)$ and $f :[0,\infty) \longrightarrow \mathbb{R}^+$  is convex with $f(0)=0$, then for all $t \geq  0$,
$$  S_{\mathcal{M}}(t)f(u) \geq f \left( S_{\mathcal{M}}(t)u \right).$$
Here, $S_{\mathcal{M}}(t)u$ is given by \eqref{conv-semigroup}, and $f(u)$ denotes the function given by $f(u)(x) = f(u(x))$.
\end{lemma}
\begin{proof} 
Let $0< s := \int_{\mathcal M} H_t(x,y)\,d\mu(y)\le 1
 $ (see estimate \eqref{energ-kernel}). Using the expression \eqref{conv-semigroup} and Jensen's inequality (see \cite[Theorem 3.3]{rudin}) with measure $d\nu = \frac{H_t(x,y)}{s}\,d\mu$, we have
	$$
	\begin{array}{ll}
		f ([S_{\mathcal{M}}(t)u](x))&= f\left ( \displaystyle \int_{\mathcal M} H_t(x,y) u(y)\,d\mu(y) \right )\\
        \noalign{\medskip}
    &=f\left ( s \displaystyle \int_{\mathcal M}  \frac{H_t(x,y)}{s} u(y)\,d\mu(y) \right ) \\
        \noalign{\medskip}
        & = f \left (s \displaystyle \int_{\mathcal M}  u(y)\,d\nu(y) + (1 - s) \,0 \right ) \\
        \noalign{\medskip}
		&\leq s \, f \left  ( \displaystyle\int_{\mathcal M}  u(y)\,d\nu(y)\right ) + (1-s) \,f(0)\\
        \noalign{\medskip}
		& \leq  s \displaystyle \int_{\mathcal M}  f(u(y) ) \,d\nu(y) \\
        \noalign{\medskip}
        & = s \displaystyle \int_{\mathcal M}  \frac{H_t(x,y)}{s} f(u(y))\,d\mu(y) = [S_{\mathcal{M}}(t)f(u)](x).
	\end{array}
	$$
    This proves the desired result.
\end{proof}

Our main result is the following theorem.

\begin{theorem}
    \label{thm1} Let $f: [0,\infty)\to [0,\infty)$ be a continuous function,
$h:(0,\infty)\to (0,\infty)$ be such that
$h \in L^1_{\mathrm{loc}}([0,\infty))\cap C((0,\infty))$,
and let $p:\mathcal{M}\to [1,\infty)$ be a bounded measurable function satisfying $1 < p^{-} \leq p(x) \leq p^{+} < + \infty$ for almost every $x \in \mathcal{M}$. We establish:
\begin{itemize}
 \item[(i)] (Global existence) 
 Assume that there exists $m>0$ such that the function $\mathcal{F}:(0,m]\to[0,\infty)$, defined by
$
\mathcal{F}(s)=\frac{1}{s}\max\{f(s)^{p^-},\,f(s)^{p^+}\},
$
is well defined for all $s\in(0,m]$. Assume moreover that $f$ and $\mathcal{F}$ are nondecreasing functions, and that there exists $v_0\in L^r(\mathcal{M} , \mu)$, $1\le r\leq \infty$, $v_0\geq  0$, $v_0 > 0$ in a set of positive measure, with $\|S_{\mathcal{M}}(t)v_0\|_{\infty} \leq m$ for all $t>0$, satisfying the condition
\begin{equation}\label{hip-princ-global}
    \int_{0}^{\infty} h(\sigma) \mathcal{F} \big(\|S_{\mathcal{M}}(\sigma)v_0\|_{\infty} \big) d\sigma < 1.
\end{equation}
Then there exists a constant $\delta > 0$ and a global solution $u \in L^{\infty}([0,\infty), L^r(\mathcal{M}, \mu) )$ to \eqref{eq1} with initial data $u_0 = \delta v_0$.
 \item[(ii)] (Nonglobal existence) Let $u_0: \mathcal{M} \to [0, \infty]$ be measurable, $u_0\geq 0$, $u_0 > 0$ in a set of positive measure. Assume that the following conditions are satisfied:
\begin{itemize}
\item[(a)] There exists a convex nondecreasing function $g:[0,\infty)\to[0,\infty)$ such that $g(0)=0$, $g(s)>0$ for all $s>0$, and
\begin{equation}\label{TT1}
f(s)\ge g(s)\quad \text{for all } s>0.    
\end{equation}
  \item[(b)] \begin{equation}\label{F_1}
     \int_{t_0}^{\infty} \frac{d \sigma}{\min \big\{g(\sigma)^{p^-}, g(\sigma)^{p^+}\big\} } < \infty
\end{equation}
for some $t_0 \geq 0$.
  \item[(c)] There exist $ \tau > 0$ such that
  \begin{equation} \label{TT2}
  \int_{ \| S_{\mathcal{M}}(\tau)u_0 \|_{\infty} }^{\infty} \frac{d \sigma}{\min \big\{g(\sigma)^{p^-}, g(\sigma)^{p^+}\big\}} < 2^{-p^+} \int_{0}^{\tau} h(\sigma) d \sigma.
  \end{equation}
\end{itemize}
Then there is no nonnegative measurable global solution $u: \mathcal{M} \times [0, \infty) \to [0, \infty]$  to the integral equation \eqref{eq2}, such that $u(x, t) < \infty $ for a.e. $(x, t) \in \mathcal{M} \times (0, \infty)$. 
\end{itemize}
\end{theorem}

\begin{remark}\label{remarks-princteo}
\begin{enumerate}
    \item[(a)] If $f$ is locally Lipschitz on $[0,\infty)$ and satisfies $f(0)=0$, then for any $m>0$ and any $0<s\le m$ we have
$$\vert f(s)-f(0) \vert \le L\vert s\vert ,
$$
where $L$ denotes a Lipschitz constant of $f$ on the interval $[0,m]$. Therefore, for all $s\in(0,m]$,
\begin{equation}\label{use-lips-localex}
    \vert f(s)^{p^{\pm}}-f(0)^{p^{\pm}}\vert 
\le p^{\pm}(mL)^{p^{\pm}-1}\vert f(s)-f(0)\vert 
\le p^{\pm}(mL)^{p^{\pm}-1}L\vert s\vert .
\end{equation}
This shows that the assumption that $f$ is locally Lipschitz with $f(0)=0$ is sufficient
to ensure that $\mathcal{F}$ is well defined on any interval $(0,m]$. Moreover, if we
assume in addition that $f$ is convex, then $\mathcal{F}$ is nondecreasing.
  
    \item[(b)] The integral condition on the reaction part in \eqref{F_1} is an example of 
a so-called Osgood-type condition, which is closely related to the 
existence of blow-up solutions (see, e.g., \cite{Laister} and \cite{meier2}).
    \item[(c)] We point out that, in the proof of Theorem \ref{thm1} (ii), we use the identity
    \begin{equation}\label{eq-sem11}
S_{\mathcal{M}}(s - t)\left( \int_0^t S_{\mathcal{M}}(t - \sigma)\varphi(\sigma,\cdot)\, d\sigma \right)
= \int_0^t S_{\mathcal{M}}(s - \sigma)\varphi(\sigma,\cdot)\, d\sigma,
\end{equation}
for $0<t<s$ and where
\begin{equation*}
    \varphi(\sigma, y) := h(\sigma)f(u(\sigma, y))^{p(y)}.
\end{equation*}
Indeed, this identity can be derived thanks to the representation of the semigroup in terms of the sub-Riemannian heat kernel $H_t:  {\mathcal{M}} \times {\mathcal{M}} \to (0, \infty)$, $t>0$, as follows. Let $u\in L^{\infty }\left ( \left ( 0,T \right ),L^{r}\left ( \mathcal{M} , \mu \right ) \right )$ be a solution of the problem (\ref{eq1}). Then the right-hand side of \eqref{eq2} belongs to the same space, and therefore we may use the kernel representation
\begin{align*}
    &\left(S_{\mathcal{M}}(s - t) \int_0^t S_{\mathcal{M}}(t - \sigma)\varphi(\sigma,\cdot)\, d\sigma \right)(x)\\
    \noalign{\medskip}
&= \int_{\mathcal{M}} H_{s - t}(x,y) \left( \int_0^t (S_{\mathcal{M}}(t - \sigma)\varphi(\sigma,\cdot))(y)\, d\sigma \right) d\mu(y).
\end{align*}
Since all functions in the integrals are measurable and nonnegative, we may apply Tonelli’s theorem \cite[Theorem 2.37]{FOLLAND},
\begin{align*}
 &\int_{\mathcal{M}} H_{s - t}(x,y)   \left( \int_0^t (S_{\mathcal{M}}(t - \sigma)\varphi(\sigma,\cdot))(y)\, d\sigma \right) d\mu(y)\\
 \noalign{\medskip}
&= \int_0^t \left( \int_{\mathcal{M}} H_{s - t}(x,y)  \, (S_{\mathcal{M}}(t - \sigma)\varphi(\sigma,\cdot))(y)\, d\mu(y) \right) d\sigma.
\end{align*}
For each fixed $\sigma$, we write
\begin{equation*}
(S_{\mathcal{M}}(t - \sigma)\varphi(\sigma,\cdot))(y) = \int_{\mathcal{M}} H_{t - \sigma}(y,z)  \, \varphi(\sigma,z)\, d\mu(z).
\end{equation*}
Substituting into the expression and using Tonelli’s theorem again, we obtain:
\begin{align*}
&\int_0^t \left( \int_{\mathcal{M}} H_{s - t}(x,y) \, (S_{\mathcal{M}}(t - \sigma)\varphi(\sigma,\cdot))(y)\, d\mu(y) \right) d\sigma \\
\noalign{\medskip}
&= \int_0^t \int_{\mathcal{M}} \left( \int_{\mathcal{M}} H_{s - t}(x,y) \, H_{t - \sigma}(y,z) \, d\mu(y) \right) \varphi(\sigma,z)\, d\mu(z)\, d\sigma.
\end{align*}
By the semigroup property of the heat kernel (Chapman-Kolmogorov equation \cite[Theorem 1.24]{Baudoin}):
\begin{equation*}
\int_{\mathcal{M}}\, H_{s - t}(x,y) \, H_{t - \sigma}(y,z) \, d\mu(y) = H_{s - \sigma}(x,z).
\end{equation*}
This establishes \eqref{eq-sem11}.
\end{enumerate}
\end{remark}

Concerning local solutions, we have the following result.
\begin{proposition}\label{LOCALEXIST}  Let $h \in L^1_{\mathrm{loc}}([0,\infty))$ be nonnegative, $f \in C([0,\infty))$ nondecreasing and locally Lipschitz continuous, with $f(0)=0$ and let $p:\mathcal{M}\to [1,\infty)$ be a bounded measurable function satisfying \eqref{ExponentVar}. Then for every $u_0 \in L^r(\mathcal{M}, \mu) $, $1 \leq r \leq \infty$, $u_0 \neq 0$, $u_0 \geq 0$, with $\|S_{\mathcal{M}}(t)u_0\|_{\infty} < \infty$ for all $t>0$, there exists a 
solution $u$ of problem \eqref{eq1} defined on $[0,T)$, for some $0<T$.
\end{proposition}
The proof of Proposition \ref{LOCALEXIST} is a straightforward adaptation of the proof of Theorem \ref{thm1} (i) and will be carried out in the next section.

\begin{remark} 
In the case $T = \infty$, we have $u\in L^{\infty }\left ( \left ( 0, \infty \right ),   L^{r}\left ( \mathcal{M} , \mu \right )  \right )$, which implies that $u(x, t) < \infty $ for almost every $(x, t) \in \mathcal{M} \times (0, \infty)$. Therefore, Theorem \ref{thm1} (ii) can be viewed as a criterion for the nonexistence of global solutions to \eqref{eq1} (in the sense of Definition \ref{Def}).
\end{remark}

\section{Proofs} \label{proooofs}

\noindent \textbf{Proof  of Theorem \ref{thm1}} (i) We employ an approach analogous to the one presented in \cite{Weissler}. Assume that the hypothesis \eqref{hip-princ-global} holds. Let $\delta > 0$ be chosen such that  
\begin{equation}\label{Con.11}
\delta = \frac{1}{\beta + 1},
\end{equation}  
where $\beta > 0$ satisfies  
\begin{equation}\label{hip-sec-global}
\int_{0}^{\infty} h(\sigma) \mathcal{F} (\|S_{\mathcal{M}}(\sigma)v_0\|_{\infty}) d\sigma < \frac{\beta}{\beta+1}<1,
\end{equation}  
for some $v_0 \in L^r(\mathcal{M}, \mu)  $ satisfying the assumptions of the theorem. Set $u_0 = \delta v_0 \in L^r(\mathcal{M}, \mu)$ and define the sequence $\{ u^k \}_{k\geq0}$ by $u^0(t) = S_{\mathcal{M}}(t)u_0$ and  
\begin{equation*}
u^k(t) = S_{\mathcal{M}}(t)u_0 + \int_{0}^{t} S_{\mathcal{M}}(t-\sigma) h(\sigma)[f(u^{k-1}(\sigma))]^{p(\cdot)}~ d \sigma,
\end{equation*}  
for $k \in \mathbb{N}$ and $t\geq 0.$

We claim that  
\begin{equation} \label{C1}
 u^{k}(t) \leq (1+\beta) S_{\mathcal{M}}(t)u_0,
\end{equation}  
for $k\geq0$ and $t>0$. To demonstrate this, we proceed by induction on $k$. Estimate (\ref{C1}) is straightforward for $k=0$. Assume that (\ref{C1}) holds for some $k \in \mathbb{N}$. 
From (\ref{Con.11}) and the assumption on $S_{\mathcal{M}}(t) v_0$, we have
$\| (1+\beta)S_{\mathcal{M}}(t)u_0 \|_{\infty}  = \| S_{\mathcal{M}}(t)v_0 \|_{\infty}   \leq  m$, for every $t>0$.
Given that $\mathcal{F}: (0,m] \to [0,\infty)$ and $f$ are nondecreasing functions, and that $s\mathcal{F}(s)=\max\{f(s)^{p^-}, f(s)^{p^+}\}$ for $s\in (0,m]$, we obtain  
\begin{align}
\begin{aligned}\label{glob-pass}
    & \quad u^{k+1}(t) \\
    \noalign{\medskip}
    &=  S_{\mathcal{M}}(t)u_0 + \int_{0}^{t} S_{\mathcal{M}}(t - \sigma) h(\sigma) \big[f(u^k(\sigma))\big]^{p(\cdot)}~d\sigma \\
     \noalign{\medskip}
 & \leq  S_{\mathcal{M}}(t)u_0 + \int_{0}^{t} h(\sigma) S_{\mathcal{M}}(t - \sigma)  \big[f((1+\beta) S_{\mathcal{M}}(\sigma)u_0)\big]^{p(\cdot)}~d\sigma \\
 \noalign{\medskip}
 & =  S_{\mathcal{M}}(t)u_0 + \int_{0}^{t} h(\sigma) S_{\mathcal{M}}(t - \sigma)  \big[f(S_{\mathcal{M}}(\sigma)v_0)\big]^{p(\cdot)}~d\sigma 
\\
 \noalign{\medskip}
 &\leq  S_{\mathcal{M}}(t)u_0 + \int_0^t h(\sigma) S_{\mathcal{M}}(t-\sigma) \max\Big\{ \big[f(S_{\mathcal{M}}(\sigma)v_0)\big]^{p^-}, \big[f(S_{\mathcal{M}}(\sigma)v_0)\big]^{p^{+}}\Big\} d\sigma\\
 \noalign{\medskip}
& =  S_{\mathcal{M}}(t)u_0 + \int_{0}^{t} h(\sigma) S_{\mathcal{M}}(t - \sigma) \mathcal{F}(S_{\mathcal{M}}(\sigma) v_0)  S_{\mathcal{M}}(\sigma) v_0d\sigma \\
\noalign{\medskip}
& \leq     S_{\mathcal{M}}(t)u_0 +  S_{\mathcal{M}}(t)v_0   \int_{0}^{t} h(\sigma)  \mathcal{F}\big(\|S_{\mathcal{M}}(\sigma) v_0\|_{\infty}\big)  ~d\sigma \\
\noalign{\medskip}
&\leq  S_{\mathcal{M}}(t)u_0 + (1+\beta) S_{\mathcal{M}}(t)u_0  ~\left( \frac{\beta}{\beta + 1} \right)= (1 + \beta) S_{\mathcal{M}}(t)u_0.
 \end{aligned}
\end{align}
Thus, claim (\ref{C1}) holds for $k+1$.

Applying induction on $k$ once more, we establish that $u^{k-1} \leq u^k$ for all $k \in \mathbb{N}$. Consequently, for every $(x,t) \in \mathcal{M} \times [0,\infty)$, the limit $\displaystyle \lim_{k \to \infty} u^k(x,t)$ exists. Let us denote this pointwise limit by $u(x,t)$. We show that the function $u: \mathcal{M} \times [0,\infty) \to (0, \infty)$, defined by $u(x,t) = \displaystyle \lim_{k \to \infty} u^k(x,t)$ is a global solution to \eqref{eq1}.

Indeed, estimate \eqref{C1}, combined with the nonnegativity of $u^k$, implies that
$u^k(t) \in L^r(\mathcal{M},\mu)$
for every $t>0$, since it is dominated by $(1+\beta) S_{\mathcal{M}}(t)u_0$ (recall from the preliminaries that the semigroup $\{S_{\mathcal{M}}(t)\}_{t \geq 0}$ is contractive on every $L^r(\mathcal{M}, \mu)$, $1\leq r \leq \infty$). Therefore, by the Dominated Convergence Theorem and the uniqueness of the limit, the function $u(\cdot, t)$ belongs to $L^r(\mathcal{M}, \mu)$ for all $t > 0$. Furthermore, since $\| (1+\beta)S_{\mathcal{M}}(t)u_0 \|_{\infty} \leq m$ for all $t > 0$, it follows that $u \in L^{\infty}([0,\infty), L^r(\mathcal{M}, \mu) )$. It remains to show that the integral equation in \eqref{eq2} is satisfied a.e. in $\mathcal{M} \times \left ( 0,\infty \right )$. 

Consider $(x, t)\in \mathcal{M} \times \left ( 0,\infty \right )$ and the sequence of functions
$$ \sigma \in [ 0,t ] \mapsto S_{\mathcal{M}}(t - \sigma) h(\sigma) \big[f(u^k(\sigma))\big]^{p(\cdot)}(x) ,\,\,\,\mbox{for} \,\,\, k \in \mathbb{N},$$
which is nonincreasing with respect to $k$. Observe that by the same estimates used in \eqref{glob-pass}, we have the bound
\begin{align}\label{dominated1}
     S_{\mathcal{M}}(t - \sigma) h(\sigma) \big[f(u^k(\sigma))\big]^{p(\cdot)}\,(x) &\leq \big[S_{\mathcal{M}}(t)v_0 \big](x)  \,  h(\sigma)  \mathcal{F}\big(\|S_{\mathcal{M}}(\sigma) v_0\|_{\infty}\big),
\end{align}
 for all $k \in \mathbb{N}$. Since the right-hand side belongs to $L^1([0,\infty))$ by \eqref{hip-sec-global}, it follows that, for every $k$, the function $$\sigma \mapsto S_{\mathcal{M}}(t - \sigma) h(\sigma) \big[f(u^k(\sigma))\big]^{p(\cdot)}(x)$$ belongs to $L^1([0,t])$. We claim that, for $\sigma \in [0,t]$, we have
\begin{equation}\label{dominated2}
    \lim_{k \to \infty} S_{\mathcal{M}}(t - \sigma) h(\sigma) \big[f(u^k(\sigma))\big]^{p(\cdot)}(x) = S_{\mathcal{M}}(t - \sigma) h(\sigma) \big[f(u(\sigma))\big]^{p(\cdot)}(x).
\end{equation}
Indeed, it follows from the continuity of $f$ that
$$\lim_{k \to \infty} H_{t -\sigma}(x,y) \,\left[f(u^{k}(y,\sigma)) \right]^{p(y)} = H_{t -\sigma}(x,y) \,\left[f(u(y,\sigma)) \right]^{p(y)}$$ for almost every $y \in \mathcal{M}$. 
Moreover, since 
$f$ is nondecreasing, we estimate similarly to \eqref{glob-pass} and obtain
\begin{align*}
    &H_{t -\sigma}(x,y) \,\left[f(u^k(y,\sigma)) \right]^{p(y)}\\
    \noalign{\medskip}
    &\leq  \max \Big\{ \big[f(\|S_{\mathcal{M}}(\sigma)v_0\|_{\infty})\big]^{p^-}, \big[f(\|S_{\mathcal{M}}(\sigma)v_0\|_{\infty})\big]^{p^{+}} \Big\} H_{t -\sigma}(x,y) .
\end{align*}
Since the right-hand side belongs to $L^1(\mathcal{M}, \mu)$, we may invoke the Dominated Convergence Theorem to conclude that
\begin{align*}
    &\lim_{k \to \infty} S_{\mathcal{M}}(t - \sigma) h(\sigma) \big[f(u^k(\sigma))\big]^{p(\cdot)}(x)\\
    \noalign{\medskip}
    &=  \lim_{k \to \infty}\int_{\mathcal{M} } H_{t -\sigma}(x,y) h(\sigma) \,\left[f(u^{k}(y,\sigma)) \right]^{p(y)}  d \mu(y) \\
    \noalign{\medskip}
    & = \int_{\mathcal{M} }\lim_{k \to \infty} H_{t -\sigma}(x,y) h(\sigma) \,\left[f(u^{k}(y,\sigma)) \right]^{p(y)}  d \mu(y) \\
    \noalign{\medskip}
    &=  S_{\mathcal{M}}(t - \sigma) h(\sigma) \big[f(u(\sigma))\big]^{p(\cdot)}(x).
\end{align*}
This proves the claim. 

Finally, in view of \eqref{dominated1} and \eqref{dominated2}, we apply the Monotone Convergence Theorem and obtain
\begin{align*}
   u(x,t) &= \lim_{k \to \infty} u^k(x,t)\\
   \noalign{\medskip}
     &  = [S_{\mathcal{M}}(t)u_0](x) + \lim_{k \to \infty} \int_{0}^{t} S_{\mathcal{M}}(t-\sigma) h(\sigma) \left[f(u^{k-1}(\sigma)) \right]^{p(\cdot)} (x)~ d \sigma\\
     \noalign{\medskip}
     &  = [S_{\mathcal{M}}(t)u_0](x) + \int_{0}^{t} S_{\mathcal{M}}(t-\sigma) h(\sigma)\left[f(u(\sigma)) \right]^{p(\cdot)} (x)~ d \sigma.
\end{align*}
This completes the proof of (i).

\medskip 
\noindent (ii) We proceed by contradiction, assuming that there exists a global solution to the integral equation
\begin{equation*}
  u(t) = S_{\mathcal{M}}(t)u_0 + \int_{0}^{t}S_{\mathcal{M}}(t - \sigma) h(\sigma)  \big[f(u(\sigma))\big]^{p(\cdot)} ~ d \sigma,  
\end{equation*}
such that $u(x, t) < \infty $ for almost every $(x, t) \in \mathcal{M} \times (0, \infty)$. Since the heat kernel $H_t$ is positive (see Section \ref{Prel}) and $f(s) \ge g(s)$ for all $s>0$, by \eqref{TT1}, we obtain
\begin{align*}
  u(t) &= S_{\mathcal{M}}(t)u_0 + \int_{0}^{t}S_{\mathcal{M}}(t - \sigma) h(\sigma)  \big[f(u(\sigma))\big]^{p(\cdot)} ~ d \sigma\\
  \noalign{\medskip}
  & \geq S_{\mathcal{M}}(t)u_0 + \int_{0}^{t}S_{\mathcal{M}}(t - \sigma) h(\sigma)  \big[g(u(\sigma))\big]^{p(\cdot)} ~ d \sigma,
\end{align*}
for almost every $(x, t) \in \mathcal{M} \times (0, \infty)$.
By the hypothesis that $u(x, s) < \infty$ for almost every $x \in \mathcal{M}$, fixing such an $s > 0$ and considering $x \in \{x \in \mathcal{M} : u(x, s) < \infty \}$, we obtain from \eqref{eq-sem11} that for any $t \in [0, s]$
\begin{equation}\label{Non1}
S_{\mathcal{M}}(s-t)u(t) \geq S_{\mathcal{M}}(s)u_0 + \int_{0}^{t} h(\sigma)S_{\mathcal{M}}(s- \sigma) \big[g(u(\sigma))\big]^{p(\cdot)}~d\sigma.
\end{equation}
For a fixed $x \in \{x \in \mathcal{M} : u(x, s) < \infty \}$, set
$$\Phi(t):=\big[S_{\mathcal{M}}(s)u_0\big](x) + \int_{0}^{t} h(\sigma)\Big[S_{\mathcal{M}}(s- \sigma) \big[g(u(\sigma))\big]^{p(\cdot)}\Big](x)~d\sigma$$
for $t\in [0,s]$. We have that $\Phi(t) \leq u(x,s) < \infty$, and the function $\Phi(t)$ is absolutely continuous on $[0, s]$. Thus, the derivative of $\Phi$ satisfies:
\begin{equation*}
  \Phi'(t) = h(t) \Big[S_{\mathcal{M}}(s-t) \big[g(u(t))\big]^{p(\cdot)}\Big](x).  
\end{equation*}
Using Lemma \ref{Ferreira-l2}, we obtain:
\begin{align*}
&\Big[S_{\mathcal{M}}(s-t) \big[g(u(t))\big]^{p(\cdot)}\Big](x)\\
\noalign{\medskip}
&= \int_{\mathcal{M}}  \big[g(u(t,y))\big]^{p(y)} \underbrace{ H_{s-t}(x,y) d \mu (y) }_{d \eta(y)}\\
\noalign{\medskip}
&\geq 2^{-p^+} \times \\
&\min\left  \{ {\left[ \int_{\mathcal{M}}  \big[g(u(t,y))\big]  H_{s-t}(x,y)d \mu (y)  \right]^{p^-}}, {\left[ \int_{\mathcal{M}}  \big[g(u(t,y))\big]  H_{s-t}(x,y) d \mu (y)  \right]^{p^+}} \right \}\\
\noalign{\medskip}
&\geq 2^{-p^+} \min\left  \{ {\big[S_{\mathcal{M}}(s-t) g(u(t)) (x)\big]^{p^-}}, {\big[S_{\mathcal{M}}(s-t) g(u(t))(x)\big]^{p^+}} \right \},   
\end{align*}
since $\textstyle \int_{\mathcal{M}} d\eta (y)=\int_{\mathcal{M}}  H_{s-t}(x,y) d \mu (y)  \leq 1 $. 
Therefore, since $p^{+} \geq p^{-} > 1$ and 
$g$ is convex nondecreasing, we can apply Lemma \ref{jensen} and \eqref{Non1} to obtain
\begin{equation}\label{Jensen}
\begin{array}{rll}
 \Phi'(t) &\geq& 2^{-p^+}\, h(t)\min\left  \{ \Big[S_{\mathcal{M}}(s-t) g(u(t))(x)\Big]^{p^-},\Big[S_{\mathcal{M}}(s-t) g(u(t))(x)\Big]^{p^+}\right \}\\
 \noalign{\medskip}
&\geq & 2^{-p^+}\, h(t) \min \left\{ \Big[g\big(S_{\mathcal{M}}(s-t)u(t)\big) (x)\Big]^{p^+}, \Big[g\big(S_{\mathcal{M}}(s-t)u(t)\big) (x)\Big]^{p^-} \right\} \\
\noalign{\medskip}
 & = & 2^{-p^+}\, h(t) \min \left\{ \big[g\big(\Phi(t)\big)\big]^{p^+}, \big[g\big(\Phi(t)\big)\big]^{p^-} \right\}.
\end{array}
\end{equation}
Define $m(t) = \min\{[g(t)]^{p^{-}}, [g(t)]^{p^{+}}\}$ for all $t \geq 0$. It is clear that $m$ is continuous on $[0, \infty)$. Then, from \eqref{Jensen}, we obtain $$\Phi'(t) \geq 2^{-p^{+}}h(t) m(\Phi(t)).$$
Next, define $$M(t) = \int_t^{+\infty} \frac{d\sigma}{m(\sigma)}.$$
Note that condition \eqref{F_1} guarantees that the function $M$ is well-defined and positive for all $t > 0$. This gives 
\begin{equation}\label{DerivativeG}
 \big[M(\Phi(t))\big]'=-\frac{\Phi'(t)}{m(\Phi(t))}\leq -2^{-p^+}\, h(t),   \end{equation}
for $0<t<s$. 
Integrating inequality \eqref{DerivativeG}  over $[0,s]$, we obtain
\begin{equation*}
  \begin{array}{ll}
-M\Big(\big[S_{\mathcal{M}}(s)u_0\big](x)\Big) &\leq M\big(\Phi(s)\big)-M\big(\Phi(0)\big)\\
\noalign{\medskip}
&=\displaystyle \int_{0}^{s}  \big[M(\Phi(\sigma))\big]'d\sigma \\
\noalign{\medskip}
& = \displaystyle - \int_{0}^{s} \frac{\Phi'(\sigma)}{m(\Phi(\sigma))} d\sigma \\
\noalign{\medskip}
&\leq\displaystyle -  2^{-p^+} \int_0^s h(\sigma) d\sigma,
\end{array}  
\end{equation*}
which is equivalent to $2^{-p^+}\int_0^s h(\sigma)d\sigma \leq M([S_{\mathcal{M}}(s)u_0](x))$. 
Observing that the left-hand side does not depend on $x$, we conclude that
$$2^{-p^+} \int_0^s h(\sigma)d\sigma \leq M\big(\|S_{\mathcal{M}}(s)u_0\|_{\infty}\big).$$
Since 
$s>0$ is arbitrary, this contradicts condition (\ref{TT2}).
\qed

\medskip
\begin{remark} Note that in the proof of Theorem \ref{thm1} (ii), $M'(t)<0$, hence $M$ is decreasing. Its inverse $M^{-1}$  is also decreasing; therefore,
    $$M\big([S_{\mathcal{M}}(s)u_0](x)\big) \geq  2^{-p^+} \int_0^s h(\sigma) d\sigma$$ 
implies
$$[S_{\mathcal{M}}(s)u_0](x) \leq M^{-1}\Big[ 2^{-p^+} \int_0^s h(\sigma) d\sigma\Big].$$
This shows, in particular, that the quantity $\|S_{\mathcal{M}}(s)u_0\|_{\infty}<\infty$ for almost every $s \in (0,T]$, $T<\infty$.
\end{remark}
\noindent \textbf{Proof of Proposition \ref{LOCALEXIST}} Following the ideas of the proof of Theorem \ref{thm1} (i), we define the sequence $\{ u^k \}_{k\geq0}$ by $u^0(t) = S_{\mathcal{M}}(t)u_0$ and  
\begin{equation*}
u^k(t) = S_{\mathcal{M}}(t)u_0 + \int_{0}^{t} S_{\mathcal{M}}(t-\sigma) h(\sigma)\big[f(u^{k-1}(\sigma))\big]^{p(\cdot)}~ d \sigma,
\end{equation*}  
for $k \in \mathbb{N}$ and $t\geq 0.$ We claim that  
\begin{equation} \label{C1-exl}
 u^{k}(t) \leq (1+ C) S_{\mathcal{M}}(t)u_0,
\end{equation}  
for all $k \geq 0$ and all $t \in (0,T)$, for some $T>0$, where $C>0$ is any fixed constant independent of $T$.  Clearly, it holds for $k=0$, and we thus assume that \eqref{C1-exl} holds for some $k$. The assumption that $f \geq 0$ is locally Lipschitz with $f(0)=0$, together with $p^+, p^- > 1$, yields estimate \eqref{use-lips-localex} of Remark \ref{remarks-princteo} for $p^+$ and $p^-$. Hence, following the computations in \eqref{glob-pass}, we obtain
\begin{align}\label{passage-local}
\begin{aligned}
     & \quad u^{k+1}(t)  \\
     \noalign{\medskip}
     &=  S_{\mathcal{M}}(t)u_0 + \int_{0}^{t} S_{\mathcal{M}}(t - \sigma) h(\sigma) \big[f(u^k(\sigma))\big]^{p(\cdot)}~d\sigma \\
     \noalign{\medskip}
 & \leq  S_{\mathcal{M}}(t)u_0 + \int_{0}^{t} h(\sigma) S_{\mathcal{M}}(t - \sigma)  \big[f((1+C) S_{\mathcal{M}}(\sigma)u_0)\big]^{p(\cdot)}~d\sigma \\
 \noalign{\medskip}
 &\leq  S_{\mathcal{M}}(t)u_0 + \int_0^t h(\sigma) S_{\mathcal{M}}(t-\sigma) \max\Big\{ \big[f(S_{\mathcal{M}}(\sigma)v_0)\big]^{p^-}, \big[f(S_{\mathcal{M}}(\sigma)v_0)\big]^{p^{+}}\Big\} d\sigma\\
 \noalign{\medskip}
& \leq  S_{\mathcal{M}}(t)u_0 + C_1 \int_{0}^{t} h(\sigma) S_{\mathcal{M}}(t - \sigma)  S_{\mathcal{M}}(\sigma) v_0d\sigma \\
\noalign{\medskip}
& \leq     S_{\mathcal{M}}(t)u_0 +  S_{\mathcal{M}}(t)v_0 \, \, C_1 \int_{0}^{t} h(\sigma) ~d\sigma \\
\noalign{\medskip}
& \leq  S_{\mathcal{M}}(t)u_0 + C_1 (1+C) \, S_{\mathcal{M}}(t)u_0 \int_{0}^{t} h(\sigma) ~d\sigma, 
\end{aligned}
\end{align}
where $v_0 = (1+C)u_0$, and $C_1 = \max \{ p^+(mL)^{p^+-1}L,\, p^-(mL)^{p^--1}L  \}$ (denotes the maximum of the constants appearing in \eqref{use-lips-localex} corresponding to $p^+$ and $p^-$). Note that we used the fact that there exists $m>0$ such that $S_{\mathcal{M}}(t)v_0 < m$, in view of the assumption that $\|S_{\mathcal{M}}(t)u_0\|_{\infty}<\infty$.  

Since $h \in L^1_{\mathrm{loc}}([0,\infty))$, we have
$
\textstyle \int_0^t h(\sigma)\,d\sigma \to 0$ as $t \to 0.$
Hence, we may choose $T>0$ such that
$\textstyle
C_1(1+C)\int_0^t h(\sigma)\,d\sigma < C $, for all  $t \in (0,T).
$
Therefore, \eqref{C1-exl} follows from this estimate together with \eqref{passage-local}.

Since \eqref{C1-exl} holds, the verification that the sequence $\{u^k\}_{k\geq 0}$ converges to a solution
$
u \in L^{\infty}\bigl((0,T),\, L^r(\mathcal{M},\mu) \bigr)
$
is similar to that in the proof of Theorem \eqref{thm1} (i), and will therefore be omitted.
\qed

\section{Applications in connected unimodular Lie groups}\label{connlieg}

In this section, we present some applications of Theorem \ref{thm1} on connected unimodular Lie groups, which provide nontrivial examples of sub-Riemannian manifolds satisfying the assumptions of our framework. Naturally, our results also apply to the Euclidean setting and, more generally, to Riemannian manifolds endowed with their natural volume measure.

Let $(\mathcal{G}, \, \circ)$ be a connected unimodular Lie group, endowed with its Haar measure $dy$, and $X = \{X_1, \dots, X_k\}$ a H\"ormander system of left-invariant vector fields. Associated with $X$, we have the Carnot-Carathéodory distance $(x, y) \in \mathcal{G} \times \mathcal{G} \to d_C(x, y)$ and denote by $d_C(x)$ the distance from the origin (the unit element of $\mathcal{G}$) to $x \in \mathcal{G}$. By \cite[Proposition III.4.2]{VA}, the behavior of the distance $d_C$ for sufficiently distant points essentially depends on $\mathcal{G}$ and not on $X$. Equipped with this H\"ormander system of left-invariant vector fields, $(\mathcal G,\circ)$ carries a complete left-invariant sub-Riemannian structure, with Haar measure as the associated smooth volume; in particular, it is a complete sub-Riemannian manifold endowed with a smooth $\sigma$-finite volume measure $dy$.

The volume function $V(t)$ represents the volume of the balls $B_{\mathcal{G}}(x, t)$ centered at $x \in \mathcal{G}$ and of radius $t > 0$ with respect to the Carnot-Carathéodory distance. It is well-known (see \cite[IV.5.8 Theorem]{VA}, \cite[Section II.4]{DUN}) that there exists an integer $d \geq 0$ (called the local dimension of $(\mathcal{G},\, X)$) and constants $c, C > 0$ such that
\begin{equation*}
   ct^d\leq  V(t) \leq C t^d, \mbox{ for } 0<t\leq 1.
\end{equation*}
The behavior of the function $t \mapsto V(t)$ for  $t \geq  1$ in the context of Lie groups can only be of two types: polynomial volume growth and exponential volume growth, independently of the choice of $X$.  We say that $\mathcal{G}$ is a group of polynomial growth if there exists an integer $D \geq 0$ (referred to as the dimension at infinity of $\mathcal{G}$) and real numbers $C, c > 0$ such that 
\begin{equation}\label{pol-grow}
   ct^D\leq  V(t) \leq C t^D
\end{equation}
uniformly for all $t \geq  1$.
We say that $\mathcal{G}$ is a group of exponential growth if there exist $\lambda, \mu > 0$ and $C, c > 0$ such that
\begin{equation*}
   ce^{\lambda t}\leq  V(t) \leq C e^{\mu t}, \quad \mbox{for all } t \geq  1.
\end{equation*}
Let $\{S_{ \mathcal{G} }(t)\}_{t \geq 0}$ be the symmetric sub-Markovian semigroup associated with the Sub-Laplacian $\Delta_{{ \mathcal{G} }}$, and let $H_t$ be its right convolution kernel such that
\begin{equation*}
    [S_{ \mathcal{G} }(t) u](x) = \int_{\mathcal{G}} H_t(y^{-1} \circ x)  u (y) dy.
\end{equation*}
From \cite[VIII Theorem 2.9]{VA}, there exist positive constants $C_1 > 0$ and $C_2 > 0$ such that  $H_t$ satisfies
\begin{equation}\label{est-ker1}
C_1 V (\sqrt{t})^{-1} \exp\left(-C_2 \frac{(d_C(x))^2}{t}\right) \leq H_t(x) \leq C_2 V (\sqrt{t})^{-1} \exp\left(-C_1 \frac{(d_C(x))^2}{t}\right)
\end{equation}
for all $t > 0$ and $x \in \mathcal{G}$, when $\mathcal{G}$ has polynomial volume growth ($H_t(x)=H_t(e\circ x)$, where $e$ denotes the identity element of $\mathcal G$). In the case where $\mathcal{G}$ has exponential volume growth, if the local dimension of $(\mathcal{G}, X)$ is $d$, then for every $n \geq d$ and $\varepsilon > 0$, there exists $C_{n, \varepsilon} > 0$ such that
\begin{equation}\label{est-ker2}
\vert H_t(x) \vert \leq C_{n, \varepsilon} t^{-\frac{n}{2}}\exp\left(-\frac{(d_C(x))^2}{(4+\varepsilon)t}\right)
\end{equation}
for all $t > 0$ and $x \in \mathcal{G}$ (see \cite[VIII Theorem 4.3]{VA}).

For connected unimodular Lie groups with polynomial volume growth, we have the following result for the case $f(u)=u$.

\begin{theorem}\label{apli-lie} Let $\mathcal{G}$ be a connected unimodular Lie group with polynomial volume growth, and suppose that $h(t) = t^q$, with $q> -1$, and $f: [0,\infty)\to [0,\infty)$ is given by $f(u) = u$.
\begin{itemize}
\item[(i)]  Let $u_0\geq 0$ be measurable and $u_0 \neq 0$. If $p^+<1+2(q+1)/D$, where $D$ is the dimension at infinity of $\mathcal{G}$, 
then problem \eqref{eq1} with initial condition $u_0$ admits no global solution.  
\item[(ii)] If $p^{-}>1+2(q+1)/D$, then problem \eqref{eq1} possesses nontrivial global solutions.
\end{itemize}
\end{theorem}

\begin{proof}(i)  
Note that $f(s)=s$ is a convex nondecreasing function with $f(0)=0$ and $f(s)>0$ for $s>0$. Since condition \eqref{F_1} holds for $f$ (for instance with $t_0=1$), provided that $1< p^- \le p(x) \le p^+$, Theorem~\ref{thm1} (ii) is therefore applicable, provided that condition~\eqref{TT2} holds for $f$ (in this case, we set $g=f$). To this end, it is sufficient to show that there exists $\tau > 0$ such that 
  \begin{equation}\label{cond-verf}
 \left(\int_{0}^{\tau} \sigma^q d \sigma \right)^{-1} \int_{ \| S_{\mathcal{\mathcal{G}}}(\tau)u_0 \|_{\infty} }^{\infty} \frac{d \sigma}{\min \{\sigma^{p^-}, \sigma^{p^+}\}} < 2^{-p^+} .
  \end{equation}

In order to estimate \eqref{cond-verf}, 
it is necessary to obtain a large-time estimate for $\| S_{\mathcal G}(t)u_0 \|_{\infty}$. Let $ u_0 > 0$ on a ball $B_{\mathcal{G}}(x_0, \delta) \subset \mathcal{G}$, where $x_0 \in \mathcal{G}$ and $\delta>0$. Then, from \eqref{est-ker1}, \eqref{pol-grow} and Fatou’s Lemma, we have
\begin{eqnarray*}
\liminf_{t \rightarrow \infty} t^{\frac{D}{2}}[ S_{\mathcal{G}}(t)u_0](x) & \geq & c\,C_1  \liminf_{t \rightarrow \infty}  \int_{\mathcal{G}} \exp \left(  - \dfrac{C_2  \left( d_C(y^{-1} \circ x) \right)^2}{t} \right) u_0(y) \, dy\\
\noalign{\medskip}
& \geq &  c\,C_1  \liminf_{t \rightarrow \infty}  \int_{B_{\mathcal{G}}(x_0, \delta)} \exp \left(  - \dfrac{C_2  \left( d_C(y^{-1} \circ x) \right)^2}{t} \right) u_0(y) \, dy\\
\noalign{\medskip}
& \geq &  c\,C_1\int_{B_{\mathcal{G}}(x_0, \delta)} u_0(y) \, dy,
\end{eqnarray*}
for any $x\in\mathcal G$. It follows that there exists a constant $C>0$ such that
$[S_{\mathcal G}(\tau)u_0](x) \ge C\,\tau^{-\frac{D}{2}}$
for all $\tau>1$ sufficiently large. Therefore,
\begin{equation}\label{est-semigrupo-unimod}
  \|S_{\mathcal{G}}(\tau)u_0 \|_{\infty} \geq C \tau^{-\frac{D}{2}}, \mbox{ for all } \tau >1  \mbox{ sufficiently large}.  
\end{equation}
Then, if $\tau>1$ is sufficiently large so that $C\,\tau^{-\frac{D}{2}}<1$, we have
\begin{align*}
    \left(\int_{0}^{\tau} \sigma^q d \sigma \right)^{-1} &\int_{ \| S_{\mathcal{\mathcal{G}}}(\tau)u_0 \|_{\infty} }^{\infty} \frac{d \sigma}{\min \{\sigma^{p^-}, \sigma^{p^+}\}}\\
    \noalign{\medskip}
    &\leq  \left(\int_{0}^{\tau} \sigma^q d \sigma \right)^{-1} \int_{ C \tau^{-D/2} }^{\infty} \frac{d \sigma}{\min \{\sigma^{p^-}, \sigma^{p^+}\}}\\
    \noalign{\medskip}
    & = \frac{\tau^{-(q+1)}}{(q+1)^{-1}} \left(  \int_{ C \tau^{-D/2} }^{1} \frac{d \sigma}{\sigma^{p^+}} + \int_{ 1}^{\infty} \frac{d \sigma}{ \sigma^{p^-}}\right)\\
    \noalign{\medskip}
    & = \frac{\tau^{-(q+1)}}{(q+1)^{-1}} \left(  \frac{C^{-p^+ +1} \tau^{\frac{D}{2}(p^+ -1)}}{p^+ -1}  - \frac{1}{p^+ -1} + \frac{1}{p^- -1}\right).
\end{align*}
Since $p^{+}<1+2(q+1)/D$, we have that  
$$ \lim_{\tau \to \infty}  \frac{\tau^{-(q+1)}}{(q+1)^{-1}} \left(  \frac{C^{-p^+ +1} \tau^{\frac{D}{2}(p^+ -1)}}{p^+ -1}  - \frac{1}{p^+ -1} + \frac{1}{p^- -1}\right) =0.$$ 
Therefore, it is possible to choose $\tau > 1$ sufficiently large (taking possibly a larger value) such that  
$$\frac{\tau^{-(q+1)}}{(q+1)^{-1}} \left(  \frac{C^{-p^+ +1} \tau^{\frac{D}{2}(p^+ -1)}}{p^+ -1}  - \frac{1}{p^+ -1} + \frac{1}{p^- -1}\right) < 2^{-p^+} .$$  
Thus, estimate \eqref{cond-verf} is verified. In view of Theorem \ref{thm1} (ii), problem \eqref{eq1} with initial condition $u_0$ does not admit a global solution.
\\
\\
(ii) Note that
$$\mathcal{F}(u)=\max \{ u^{p^- -1}, u^{p^+ -1}\}$$
is nondecreasing in $(0, \infty)$. Let $v_0\in L^{\infty}(\mathcal{G})$ 
 be any nonnegative compactly supported function whose support is contained in a ball $B_{\mathcal{G}}(x_0, R)$ with $x_0\in \mathcal{G}$ and $R>0$, satisfying
$$ \|v_0\|_{\infty} < \min\left\{ 1, \, \dfrac{1}{CC_2 V(R) }  \left( \dfrac{q \left[ \frac{D}{2}(p^- -1)-(q+1) \right]}{(q+1)  } \right)^{\frac{1}{p^- -1}}\right\},$$
where $V(R)>0$ denotes the volume of the ball $B_{\mathcal{G}}(x_0, R)$, $C>0$ is the constant from \eqref{pol-grow}, $C_2>0$ is the constant from \eqref{est-ker1}, and $\frac{D}{2}(p^- -1)-(q+1)>0$, since $1 + \frac{2(q+1)}{D} < p^{-}$. Then, from \eqref{est-ker1} and \eqref{pol-grow}, we have for $t>1$
\begin{align}\label{citar-pro1}
    \begin{aligned}
        S_{\mathcal{G}}(t)v_0(x) & \leq  CC_2 t^{-\frac{D}{2}}  \int_{\mathcal{G}} \exp \left(  - \dfrac{C_1  \left( d_C(y^{-1} \circ x) \right)^2}{t} \right) v_0(y) \, dy\\
        \noalign{\medskip}
& \leq   CC_2 t^{-\frac{D}{2}}  \int_{\mathcal{G}}  v_0(y) \, dy\\
\noalign{\medskip}
& \leq   CC_2 t^{-\frac{D}{2}} \|v_0\|_{\infty} V(R) \\
\noalign{\medskip}
&<  \left( \dfrac{q \left[ \frac{D}{2}(p^- -1)-(q+1) \right]}{(q+1)  } \right)^{\frac{1}{p^- -1}}  t^{-\frac{D}{2}}.
    \end{aligned}
\end{align}
From this, it follows that 
\begin{equation}\label{decay-1}
    \|S_{\mathcal{G}}(t)v_0 \|_{\infty} < C^\star t^{-\frac{D}{2}} \mbox{ for } t>1,
\end{equation} 
where $C^\star = \left( \frac{q \left[ \frac{D}{2}(p^- -1)-(q+1) \right]}{(q+1)  } \right)^{\frac{1}{p^- -1}}.$ 

Thus, by the $L^\infty$-contractivity and since $1+\frac{2(q+1)}{D}<p^{-}$, we obtain
\begin{align*}
&\int_{0}^{\infty} h(\sigma) \mathcal{F} \left(\|S_{\mathcal{G}}(\sigma)v_0\|_{\infty} \right) d\sigma
\end{align*}

\begin{align}\label{citar-prox2}
\begin{aligned}
    &= \left( \int_{0}^{1} + \int_{1}^{\infty} \right) \sigma^q  \max \{ \|S_{\mathcal{G}}(\sigma)v_0\|_{\infty}^{p^- -1}, \|S_{\mathcal{G}}(\sigma)v_0\|_{\infty}^{p^+ -1}\} d\sigma \\
    \noalign{\medskip}
    & <  \max \dfrac{\{ \|v_0\|_{\infty}^{p^- -1}, \|v_0\|_{\infty}^{p^+ -1}\}}{q+1} +  (C^\star )^{(p^- -1)}\int_{1}^{\infty}  \sigma^{ - \frac{D}{2}(p^- -1)+q} d\sigma\\
    \noalign{\medskip}
    & \leq  \dfrac{\|v_0\|_{\infty}^{p^- -1}}{q+1} +  \dfrac{(C^\star )^{(p^- -1)}}{ ({D}/{2})(p^- -1)-(q+1) }\\
    \noalign{\medskip}
    & < \dfrac{1}{q+1} + \dfrac{q}{q+1} =1.
\end{aligned}
\end{align}
It follows from Theorem \ref{thm1} (i) that there exists $\delta>0$ such that the solution of problem (\ref{eq1}) with initial datum $u_0 = \delta v_0$ is global.
\end{proof}

\begin{remark}
\begin{enumerate}
    \item[(a)] When $p(x)=p$ is constant and $q=0$, Theorem \ref{apli-lie} recovers the Fujita critical exponent stated in Theorem \ref{Rush1} (see \cite[Theorem 1.3]{RY}) for connected unimodular Lie groups with polynomial volume growth.
    \item [(b)]  Note that Theorem \ref{apli-lie} does not address the case where $p^- \leq 1+2(q+1)/D \leq p^+$. This regime is more delicate; therefore, we devote Subsection \ref{crit-case-section} to showing that, in the setting of homogeneous Carnot groups, there exists a class of initial data for which no global solution exists. In particular, our result is consistent with the critical case $p^- = p^+ = 1+2(q+1)/D$ established in \cite[Theorem 1.3 (ii)]{RY}.
\end{enumerate}

\end{remark}

Concerning connected unimodular Lie groups with exponential volume growth, we have the following result:

\begin{theorem}\label{lieee22222} Let $\mathcal{G}$ be a connected unimodular Lie group with exponential volume growth, $X$ a H\"ormander system of left-invariant vector fields, $h(t) = t^q$, $q >-1$ and $f(u) = u$. If $p^{-}>1+2(q+1)/d$, where $d$ is the local dimension of $(\mathcal{G}, X)$, then problem \eqref{eq1} possesses nontrivial global solutions.
\end{theorem}
\begin{proof} The proof is essentially the same as that of Theorem \ref{apli-lie} (ii). We consider a nonnegative compactly supported initial datum $v_0$ with $\|v_0\|_{\infty}$ sufficiently small. By using estimate \eqref{est-ker2} instead of \eqref{est-ker1} in \eqref{citar-pro1}, we obtain that 
\begin{equation*}
    \|S_{\mathcal{G}}(t)v_0\|_{\infty} < C^\star t^{-\frac{d}{2}} \mbox{ for } t>1.
\end{equation*}
 With this estimate, and using the fact that $p^{-}>1+2(q+1)/d$, we can estimate $\textstyle\int_{0}^{\infty} h(\sigma) \mathcal{F} \left(\|S_{\mathcal{G}}(\sigma)v_0\|_{\infty} \right) d\sigma$ in the same way as in \eqref{citar-prox2}. This completes the proof.
\end{proof}

The heat equation with exponential reaction, known as the Gelfand equation, was introduced by Frank-Kamenetskii~\cite{frank-kamenetskii} as a model for thermal explosions and solid fuel combustion.
Its blow-up behavior has been widely studied; see, e.g., Gelfand~\cite{gelfand} and Kaplan~\cite{kaplan}. As a consequence of Theorem~\ref{thm1} (ii), we obtain the following nonexistence result for global solutions when $f$ and $h$ are of exponential type.

\begin{theorem}\label{apli-lie-case-exp} Let $\mathcal{G}$ be a connected unimodular Lie group with polynomial volume growth, and suppose that $h(t) = e^{\theta t}$ for some $\theta>0$,  and $f(u) = e^{\gamma u} - 1$ with $\gamma>0$.  If the initial data $u_0$ is measurable, $u_0\geq 0$, $u_0 \neq 0$, then problem \eqref{eq1} with initial condition $u_0$ admits no global solution.
\end{theorem}
\begin{proof} It is clear that $f(0) = 0$ and $f(s)>0$ for $s>0$. Moreover, for all $\sigma \geq \ln(2)/\gamma$, we have
$$
(e^{\gamma \sigma} - 1)^{p^+} \geq 2^{-p^+} e^{p^+\gamma \sigma} \quad \text{and} \quad (e^{\gamma \sigma} - 1)^{p^-} \geq 2^{-p^-} e^{p^-\gamma \sigma}.
$$
Therefore, 
\begin{align*}
&\int_{\ln(2)/\gamma}^{\infty} \frac{d\sigma}{ \min\Big\{(e^{\gamma \sigma} - 1)^{p^+}, (e^{\gamma \sigma} - 1)^{p^-} \Big\}} \\
\noalign{\medskip}
    &\leq {2^{1+p^+}} \int_{\ln(2)/\gamma}^{\infty} \big[e^{-p^+\gamma \sigma} + e^{-p^-\gamma\sigma}  \big]\,d\sigma \\
    \noalign{\medskip}
    &\leq 2^{1+p^+}\left( \frac{e^{-p^+\ln2}}{\gamma p^+} + \frac{e^{-p^-\ln2}}{\gamma p^-}  \right)<\infty,
\end{align*}
and hence, condition (\ref{F_1}) holds. Furthermore, since $f$ is convex, in order to apply Theorem \ref{thm1} (ii) with $g=f$, it remains to show that condition \eqref{TT2} holds. To this end, we estimate using the semigroup estimate from \eqref{est-semigrupo-unimod}:
\begin{align*}
    \begin{aligned}
        &\left(\int_{0}^{\tau} e^{ \theta\sigma}  d \sigma \right)^{-1} \int_{ \| S_{\mathcal{\mathcal{G}}}(\tau)u_0 \|_{\infty} }^{\infty} \frac{d\sigma}{ \min\Big\{(e^{\gamma \sigma} - 1)^{p^+}, (e^{\gamma \sigma} - 1)^{p^-} \Big\}} \\
        \noalign{\medskip}
       &\leq  {2^{1+p^+}} \left(\int_{0}^{\tau} e^{ \theta\sigma}  d \sigma \right)^{-1} \int_{\| S_{\mathcal{\mathcal{G}}}(\tau)u_0 \|_{\infty}}^{\infty} \big[e^{-p^+\gamma \sigma} + e^{-p^-\gamma\sigma}  \big]\,d\sigma\\
       \noalign{\medskip}
       & \leq \left(  \dfrac{2^{1+p^+} \theta}{e^{\theta \tau } -1} \right) \int_{ C \tau^{-D/2} }^{\infty} \big[e^{-p^+\gamma \sigma} + e^{-p^-\gamma\sigma}  \big]\,d\sigma\\
       \noalign{\medskip}
       & \leq  \left(  \dfrac{2^{1+p^+} \theta}{e^{\theta \tau } -1} \right) \left( \frac{e^{-Cp^+ \gamma \tau^{-D/2}}}{\gamma p^+} + \frac{e^{-Cp^- \gamma \tau^{-D/2}}}{\gamma p^-}  \right)\\
       \noalign{\medskip}
       & =   \left(  \dfrac{2^{1+p^+} \theta}{e^{\theta \tau } -1} \right) \left( \frac{e^{-Cp^+ \gamma \tau^{-D/2}}}{\gamma p^+} + \frac{e^{-Cp^- \gamma \tau^{-D/2}}}{\gamma p^-}  \right).
    \end{aligned}
\end{align*}
  Therefore, since the limit of the last expression is zero as $\tau\to\infty$, it is possible to choose $\tau>1$ sufficiently large such that 
$$\left(  \dfrac{2^{1+p^+} \theta}{e^{\theta \tau } -1} \right) \left( \frac{e^{-Cp^+ \gamma \tau^{-D/2}}}{\gamma p^+} + \frac{e^{-Cp^- \gamma \tau^{-D/2}}}{\gamma p^-}  \right) < 2^{-p^+} .$$  
In view of Theorem \ref{thm1} (ii), there is no global solution to problem \eqref{eq1} with initial condition $u_0$.
\end{proof}

It is worth noting that the function $h(t)$ in Theorem \ref{apli-lie-case-exp} is only locally integrable. In the next theorem, we show that if $h(t)$ is integrable on $[0, \infty)$, then global solutions to problem (\ref{eq1}) can be obtained for a certain class of exponential nonlinearities.

\begin{theorem}\label{lieee222223} Let $\mathcal{G}$ be a connected unimodular Lie group with either polynomial or exponential volume growth. Assume that $h(t)$ is integrable on $[0, \infty)$ and $f(u) = u^{\alpha}e^{ u^{\beta}}$ for some $\beta > 0$ and $\alpha>1/p^-$. Then, problem \eqref{eq1} admits nontrivial global solutions.
\end{theorem}
\begin{proof} Note that $\mathcal{F}(u)$ is nondecreasing in $(0, \infty)$. Since $\|S_{\mathcal G}(t)v_0\|_{\infty}\le \|v_0\|_{\infty}$ for $v_0 \in L^{\infty}(\mathcal{G} )$, and given that
$$
\lim_{t \to 0^+} \frac{1}{t} \max \Big\{ t^{\alpha p^+}e^{ p^+ t^{\beta}}  , t^{\alpha p^-}e^{ p^- t^{\beta}}  \Big\} = 0 ,
$$
it follows that
\begin{equation*}
\dfrac{  \max \Big\{ \big\|S_{\mathcal{G}}(\sigma)v_0\big\|_{\infty}^{\alpha p^+}e^{ p^+ \big\|S_{\mathcal{G}}(\sigma)v_0\big\|_{\infty}^{\beta}}  , \big\|S_{\mathcal{G}}(\sigma)v_0\big\|_{\infty}^{\alpha p^-}e^{ p^- \big\|S_{\mathcal{G}}(\sigma)v_0\big\|_{\infty}^{\beta}}  \Big\} }{{\big\|S_{\mathcal{G}}(\sigma)v_0\big\|_{\infty}}} <  \dfrac{1}{\displaystyle \int_{0}^{\infty} h(\sigma) d\sigma +\varepsilon },
\end{equation*}
for some $\varepsilon > 0$, where $v_0$ is a nonnegative initial datum with sufficiently small $\|v_0\|_{\infty}$. Thus, 
\begin{align*}
\begin{aligned}
   & \int_{0}^{\infty} h(\sigma) \mathcal{F} \left(\|S_{\mathcal{G}}(\sigma)v_0\|_{\infty} \right) d\sigma\\
   \noalign{\medskip}
    &=  \int_{0}^{\infty}  h(\sigma)  \dfrac{  \max \Big\{ \|S_{\mathcal{G}}(\sigma)v_0\|_{\infty}^{\alpha p^+}e^{ p^+ \|S_{\mathcal{G}}(\sigma)v_0\|_{\infty}^{\beta}}  , \|S_{\mathcal{G}}(\sigma)v_0\|_{\infty}^{\alpha p^-}e^{ p^- \|S_{\mathcal{G}}(\sigma)v_0\|_{\infty}^{\beta}}  \Big\} }{{\|S_{\mathcal{G}}(\sigma)v_0\|_{\infty}}} d\sigma \\
    \noalign{\medskip}
    & <  \dfrac{\displaystyle \int_{0}^{\infty} h(\sigma) d\sigma}{ \displaystyle \int_{0}^{\infty} h(\sigma) d\sigma +\varepsilon }<1.
\end{aligned}
\end{align*}
In view of Theorem \ref{thm1} (i), there exists a constant $\delta > 0$ such that, setting $u_0 = \delta v_0$, the solution to problem (\ref{eq1}) is global.
\end{proof}

The next result addresses a case where the nonlinearity $f$ is not convex.

\begin{theorem}\label{apli-lie-case-conv} 
Let $\mathcal{G}$ be a connected unimodular Lie group with polynomial volume growth, and suppose that $h(t)=t^q$ for some $q>-1$. Let $\alpha>1$ and $0<\beta<1/p^+$, and 
$$
f(u)=
\begin{cases}
\,ku^{\alpha}, & 0\le u\le 1/2,\\
\noalign{\medskip}
\, \left(A_1u\ - B_1\right)^{\beta} & 1/2\leq u\le 1,\\
\noalign{\medskip}
\, k\alpha \left(\frac{1}{2}\right)^{\alpha -1} u - k\alpha \left(\frac{1}{2}\right)^{\alpha } + k \left(\frac{1}{2}\right)^{\alpha }, &  u\ge 1,
\end{cases}
$$
where $k>0$,
\begin{align*}
    A_1 &= 2\left[ k (\alpha+1)\left(\frac{1}{2}\right)^{\alpha}\right]^{\frac{1}{\beta}}  -2\left[ k\left(\frac{1}{2}\right)^{\alpha}\right]^{\frac{1}{\beta}},\\
    \noalign{\medskip}
     B_1 &= \left[ k (\alpha+1)\left(\frac{1}{2}\right)^{\alpha}\right]^{\frac{1}{\beta}}  -2\left[k\left(\frac{1}{2}\right)^{\alpha}\right]^{\frac{1}{\beta}}.
\end{align*}
\begin{itemize}
\item[(i)]  If $$p^{+}<\frac{1}{\alpha}+\frac{2(q+1)}{ \alpha D},$$
and the initial data $u_0$ is measurable, $u_0 \geq  0$, $u_0 \neq 0$, then problem \eqref{eq1} with initial condition $u_0$ does not admit a global solution. 
\item[(ii)] If $$p^-> \frac{1}{\alpha}+\frac{2(q+1)}{ \alpha D},$$ then problem \eqref{eq1} possesses nontrivial global solutions.
\end{itemize}

\end{theorem}
\begin{proof} (i) Consider
 $$
g(u)=
\begin{cases}
\,ku^{\alpha}, & 0\le u\le 1/2,\\
\noalign{\medskip}
\, k\alpha \left(\frac{1}{2}\right)^{\alpha -1} u - k\alpha \left(\frac{1}{2}\right)^{\alpha } + k \left(\frac{1}{2}\right)^{\alpha }, &  u\ge 1/2.
\end{cases}
$$
Since the left and right derivatives of $g$ exist at the junction point and satisfy
$$ g'_-\left(\frac{1}{2} \right) = \alpha k \left(\frac{1}{2}\right)^{\alpha -1} = g'_+\left(\frac{1}{2} \right),$$
then $g$ is convex. Moreover, note that $f(s) \geq g(s)$, $g(s)>0$ for all $s>0$, $g(0)=0$, hence \eqref{TT1} holds. Condition \eqref{F_1} is satisfied. Indeed, for any $t_0>1$ such that
$C_1\sigma+C_2>1$ for all $\sigma\in [t_0,\infty)$, we have
\begin{equation*}
     \int_{t_0}^{\infty} \frac{d \sigma}{\min \{g(\sigma)^{p^-}, g(\sigma)^{p^+}\} } = \int_{t_0}^{\infty} \frac{d \sigma}{(C_1\sigma+C_2)^{p^-}}  < \infty,
\end{equation*}
where $C_1=k\alpha \left(\frac{1}{2}\right)^{\alpha -1}$ and $C_2 = - k\alpha \left(\frac{1}{2}\right)^{\alpha } + k \left(\frac{1}{2}\right)^{\alpha }$. 

Now, we verify that condition \eqref{TT2} is fulfilled. Indeed, choosing $\tau>1$ sufficiently large in estimate \eqref{est-semigrupo-unimod} so that $C \tau^{-D/2} < \tfrac12$, we obtain
\begin{align*}
    &\left(\int_{0}^{\tau} \sigma^q d \sigma \right)^{-1}\int_{ \| S_{\mathcal{M}}(\tau)u_0 \|_{\infty} }^{\infty} \frac{d \sigma}{\min \{g(\sigma)^{p^-}, g(\sigma)^{p^+}\}} \\
    \noalign{\medskip}
    &\leq  \left(\int_{0}^{\tau} \sigma^q d \sigma \right)^{-1} \int_{ C \tau^{-D/2} }^{\infty} \frac{d \sigma}{\min \{g(\sigma)^{p^-}, g(\sigma)^{p^+}\}}\\
    \noalign{\medskip}
    & = \frac{\tau^{-(q+1)}}{(q+1)^{-1}} \times \\
    \noalign{\medskip}
 &   \left( \dfrac{1}{\min\{ k^{p^{-}} , k^{p^{+}}\}}  \int_{ C \tau^{-D/2} }^{1/2} \frac{d \sigma}{\sigma^{\alpha p^+}} +  \int_{ 1/2 }^{t_0} \frac{d \sigma}{(C_1\sigma+C_2)^{p^+}} + \int_{ t_0 }^{\infty} \frac{d \sigma}{(C_1\sigma+C_2)^{p^-}}  \right)\\
 \noalign{\medskip}
    & = \frac{\tau^{-(q+1)}}{(q+1)^{-1}} \left( C_3 \tau^{\frac{D}{2}(\alpha p^+ -1)} +C_4\right),
\end{align*}
where $C_3, \, C_4$ are constants independent of 
$\tau$. Since $\alpha p^{+}<1+2(q+1)/D$, we have that  
$$ \lim_{\tau \to \infty}   \frac{\tau^{-(q+1)}}{(q+1)^{-1}} \left( C_3 \tau^{\frac{D}{2}(\alpha p^+ -1)} +C_4\right) =0.$$ 
Therefore, it is possible to choose $\tau > 1$ (possibly larger than before) such that 
$$  \frac{\tau^{-(q+1)}}{(q+1)^{-1}} \left( C_3 \tau^{\frac{D}{2}(\alpha p^+ -1)} +C_4\right)< 2^{-p^+} .$$  
Thus, estimate \eqref{TT2} holds. By Theorem \ref{thm1} (ii), problem \eqref{eq1} with initial condition $u_0$ does not admit a global solution.
\\
~~\\
(ii) Since $f(u)^p = k u^{p \alpha}$ for $0\leq u \leq 1/2$, and $\|S_{\mathcal{G}}(t)v_0 \|_{\infty} < 1/2$ for all $t>0$ provided that $v_0$ is sufficiently small, the proof follows the same argument as that of Theorem~\ref{apli-lie}~(ii).
\end{proof}

\subsection{A result for a nonconvex nonlinearity without a polynomial lower bound}

Theorem~\ref{thm1} (ii) recovers the results of Theorem~\ref{Rush1} in the particular case $p(x)=p>1$ and $g(u)=C_0u^{q_0}$, with $q=q_0/p>1$. Moreover, it extends to a broader context. For instance, consider the nonconvex case where
\begin{equation}\label{NoRush}
f(u)=(u+1)(\ln(u+1))^{\alpha} + \sin^2(u),
\end{equation}
with $p(x)=p$ and $\alpha>1$. In this setting, there is no polynomial lower bound of the form $C_0 u^{q_0}$ with $q_0/p>1$ and $C_0>0$ such that
\begin{equation}\label{PolyRush}
f(u)\geq C_0 u^{q_0/p} \qquad \text{for all } u\geq 0,
\end{equation}
and therefore the nonglobal existence result in Theorem~\ref{Rush1} does not apply. Indeed, since $q=q_0/p>1$, we have
\begin{eqnarray*}
 & & \lim_{u \to +\infty} \frac{(u+1)(\ln(u+1))^{\alpha} + \sin^2(u)}{u^q} \\
 \noalign{\medskip}
 &=& \lim_{u \to +\infty}  \left(1 + \frac{1}{u}\right) \frac{(\ln(u+1))^{\alpha}}{u^{q-1}} + \lim_{u \to +\infty} \frac{\sin^2(u)}{u^{q-1}} =0,
\end{eqnarray*}
thus the polynomial lower bound \eqref{PolyRush} fails whenever $q_0/p>1$. 

Note also that \eqref{PolyRush} fails in the case $q=q_0/p \leq 1$. Indeed
\begin{eqnarray*}
0& \leq& \lim_{u \to 0^+} \frac{(u+1)(\ln(u+1))^{\alpha} + \sin^2(u)}{u^q}   \\
\noalign{\medskip}
& \leq  & \lim_{u \to 0^+} C \left[ u^{\alpha -q} + u^{2 -q} \right] = 0. 
\end{eqnarray*}

In the following theorem, we consider problem~\eqref{eq1} with $f$ defined by~\eqref{NoRush} and $p(x)=p>1$, and establish a Fujita-type result for the case $1<\alpha<2$.

\begin{theorem}\label{log-exam} Let $p(x)=p>1$ constant, $\alpha>1$, $\mathcal{G}$ be a connected unimodular Lie group with polynomial volume growth, $h(t) = t^q$, with $q> -1$, and $f: [0,\infty)\to [0,\infty)$ is given by $$f(u) = (u+1)(\ln(u+1))^{\alpha} + \sin^2(u)$$
\begin{itemize}
\item[(i)]  Let $u_0\geq 0$ be measurable, and $u_0 \neq 0$. If $$p <\frac{1}{\alpha}+\frac{2(q+1)}{ \alpha D},$$ where $D$ is provided in \eqref{pol-grow}, 
then problem \eqref{eq1} with initial condition $u_0$ admits no global solution.  
\item[(ii)] If $1<\alpha<2$ and
$$p> \frac{1}{\alpha}+\frac{2(q+1)}{ \alpha D},$$
then problem \eqref{eq1} possesses nontrivial global solutions.
\end{itemize}
\end{theorem}
\begin{proof}
(i) Let's consider the following function $g: [0,\infty) \rightarrow [0,\infty)$ defined by $$g(u) = (u+1)(\ln(u+1))^{\alpha}.$$
Then $g$ is a convex and nondecreasing function on $[0,\infty)$ such that $g(0) = 0$ and $0\leq g(u) \leq f(u)$ for all $u \geq 0$. Moreover 
$$\displaystyle\int_{t_0}^{\infty} \frac{du}{g(u)^p} \leq \displaystyle\int_{t_0}^{\infty} \frac{du}{(u+1)^p} < \infty \quad \mbox{for all}\,\, t_0 > e-1.$$ 
To apply Theorem~\ref{thm1}, it remains only to verify that condition \eqref{TT2} holds. Since 
$$
\lim_{t \to 0^+} \frac{\ln(1+t)}{t} = 1,
$$
there exist positive constants $\delta$ and $c_0$, sufficiently small, such that
$$
\ln(1+t) \geq c_0 \, t, \quad \text{for all } t \in [0,\delta).
$$
By the asymptotic estimate for the semigroup in \eqref{est-semigrupo-unimod}, if $\tau>1$ is chosen sufficiently large so that $C\,\tau^{-\frac{D}{2}}<\delta$, then we have
\begin{eqnarray*}
&&\left[\int_{ \| S_{\mathcal{G}}(\tau)u_0 \|_{\infty} }^{\infty} \frac{d \sigma}{g(\sigma)^{p}} \right] \left(\int_{0}^{\tau} \sigma^q d \sigma\right)^{-1} \\
\noalign{\medskip}
&=& \left[\int_{ \| S_{\mathcal{G}}(\tau)u_0 \|_{\infty} }^{\infty} \big((\sigma +1)(\ln(\sigma+1))^{\alpha}\big)^{-p}d \sigma \right] \, \left((q+1) \tau^{-(q+1)}\right)  \\
\noalign{\medskip}
&\leq& \left[\int_{ \| S_{\mathcal{G}}(\tau)u_0 \|_{\infty} }^{\delta} (\ln(\sigma+1))^{-p\alpha}d \sigma \right] \, \left((q+1) \tau^{-(q+1)}\right)
 \\
\noalign{\medskip}
&& +  \left[\int_{\delta}^{\infty} (\sigma + 1)^{-p} (\ln(\sigma+1))^{-p\alpha}d \sigma \right] \, \left((q+1) \tau^{-(q+1)}\right) \\
\noalign{\medskip}
&\leq& \left[c_{0}^{-p\alpha}\int_{C \tau^{-\frac{D}{2}}}^{\delta} \sigma^{-p\alpha}d \sigma \right] \, \left((q+1) \tau^{-(q+1)}\right) \\
\noalign{\medskip}
&& + \frac{(\delta+1)^{1-p}}
{(p\alpha-1) \,(\ln(\delta+1))^{p\alpha-1}} \, \left((q+1) \tau^{-(q+1)}\right) \\
\noalign{\medskip}
&=&  \left[\frac{c_{0}^{-p\alpha}(q+1)}{p\alpha-1}       \left( C^{1-p\alpha} \tau^{\frac{D (p\alpha-1)}{2}-(q+1)} - \delta^{1-p\alpha} \tau^{-(q+1)} \right)  \right]  \\
\noalign{\medskip}
&& + \frac{(q+1)(\delta+1)^{1-p}}
{(p\alpha-1) \,(\ln(\delta+1))^{p\alpha-1}} \,  \tau^{-(q+1)}.
\end{eqnarray*}
Since $p \alpha <1+ {2(q+1)}/{ D}$, we may increase the previously chosen value of $\tau$, if necessary, so that
$$
\left[\int_{ \| S_{\mathcal{G}}(\tau)u_0 \|_{\infty} }^{\infty} \frac{d \sigma}{g(\sigma)^{p}} \right] \left(\int_{0}^{\tau} \sigma^q d \sigma\right)^{-1}
\leq 2^{-p^{+}}.
$$
The conclusion now follows from Theorem~\ref{thm1} (ii).
\\
~~\\
(ii) Since $1<\alpha<2$, there exists a constant $C_0 > 0$ such that
$$
f(u)^p = \big[(u+1) \big(\ln(u+1)\big)^{\alpha} + \sin^2(u)\big]^p \;\leq\; C_0\,u^{p \alpha} , 
\quad \text{for } 0<u <1.
$$
Since $\|S_{\mathcal{G}}(t)v_0 \|_{\infty} < 1$ for all $t>0$, provided that $v_0$ is sufficiently small, the proof then proceeds as in Theorem~\ref{apli-lie}~(ii).
\end{proof}

\section{Applications in homogeneous Carnot groups}\label{SectionHeisen}

In this section we assume that $\mathbb{R}^N$ is equipped with a Lie group
structure given by the composition law $\circ$. Moreover, we assume that
$\mathbb{R}^N$ is endowed with a homogeneous group structure through a family
$\{\delta_\lambda\}_{\lambda>0}$ of automorphisms of $(\mathbb{R}^N,\circ)$,
called dilations, of the form
$$
\delta_\lambda(x^{(1)},x^{(2)},\ldots,x^{(r)})
=
(\lambda x^{(1)},\lambda^2 x^{(2)},\ldots,\lambda^r x^{(r)}).
$$

Here $x^{(i)}\in\mathbb{R}^{N_i}$ for $i=1,\ldots,r$, with
$$
N_1+N_2+\cdots+N_r=N.
$$

Let $\mathfrak g$ denote the Lie algebra of the group $(\mathbb{R}^N,\circ)$.
For $i=1,\ldots,N_1$, let $X_i$ be the unique vector field in $\mathfrak g$
which agrees with $\partial/\partial x_i$ at the origin.

If the Lie algebra generated by  $X_1,\ldots,X_{N_1}$
 coincides with $\mathfrak g$,
then $\mathcal{G}=(\mathbb{R}^N,\circ,\delta_\lambda)$ is called a homogeneous Carnot group.

In comparison with the previous section, the essential difference is that we now impose a stratified structure on the group. Recall that a stratified group (or Carnot group) $\mathcal{G}$ is a connected and simply connected Lie group whose Lie algebra $\mathfrak g$ admits a stratification, that is, a direct sum decomposition
$$
\mathfrak g = V_1 \oplus V_2 \oplus \cdots \oplus V_r
$$
such that
$$
[V_1,V_{i-1}] = V_i, \qquad 2 \le i \le r.
$$

In particular, stratified groups are unimodular. Moreover, it is shown in \cite[Theorem 2.2.18]{BONFIGLIOLI} that every stratified group is isomorphic to a homogeneous Carnot group. Therefore, without loss of generality, a stratified group can be identified with a homogeneous Carnot group.

In this setting, the Lebesgue measure on $\mathbb{R}^N$ coincides with the Haar measure of the group $\mathcal{G}=(\mathbb{R}^N,\circ,\delta_\lambda)$. In particular, it is invariant under group translations and homogeneous with respect to the dilations $\{\delta_\lambda\}_{\lambda>0}$, namely
$$
\vert\delta_\lambda(E) \vert=\lambda^Q \vert E\vert
$$
for every measurable set $E\subset\mathbb{R}^N$, where
$$
Q=\sum_{j=1}^{r} j \, N_j
$$
is called the homogeneous dimension of $\mathcal{G}$. For homogeneous Carnot groups, the dimension at infinity 
$D$ coincides with the homogeneous dimension 
$Q$ (see \cite[IV Definition 5.4]{VA}). In particular, the volume of the balls $B_{\mathcal{G}}(x,t)$ with respect to the Carnot--Carathéodory distance satisfies
\begin{equation}\label{pol-grow-hom}
    V(t) = C t^Q, \quad \mbox{for all} \,\,t>0.
\end{equation}

In a homogeneous Carnot group $\mathcal{G}$, a homogeneous norm is a continuous function
$\vert \cdot \vert_{\mathcal{G}}:\mathcal{G}\to[0,\infty)$
such that
$$
\vert \delta_\lambda(x) \vert_{\mathcal{G}} =\lambda\, \vert x \vert_{\mathcal{G}} , \qquad \text{for every } \lambda>0 \text{ and } x\in \mathcal{G};
$$
and
$$
\vert x \vert_{\mathcal{G}} >0 \quad \text{if and only if} \quad x\neq 0.
$$
The function
\begin{equation}\label{hom-norm-G}
    \vert x\vert _{\mathcal{G}}=\left(\sum_{j=1}^{r} \vert x^{(j)}\vert^{\frac{2r!}{j}}\right)^{\frac{1}{2r!}},
\end{equation}
is an example of a homogeneous norm on $\mathcal{G}$, smooth outside the origin, where $\vert x^{(j)} \vert$ denotes the Euclidean norm in $\mathbb{R}^{N_j}$. It is well known that all homogeneous norms on $\mathcal{G}$ are equivalent (in particular $\vert \cdot \vert_{\mathcal{G}}$ is equivalent to the Carnot--Carathéodory distance $d_C(\cdot)$) \cite[Proposition 5.1.4 and Theorem 5.2.8]{BONFIGLIOLI}. Moreover, the inverse map and the group operation satisfy pseudo-triangle inequalities, namely there exist constants $C, \,c>0$ such that $c^{-1}\, \vert x\vert _{\mathcal{G}} \, \leq \, \vert x^{-1}\vert_{\mathcal{G}} \, \leq \, c\, \vert x\vert_{\mathcal{G}}$ and $\vert  y^{-1} \circ x\vert _{\mathcal{G}} \, \leq \, C( \vert x \vert_{\mathcal{G}} + \vert y\vert_{\mathcal{G}})$ for all $x,y\in G$ \cite[Corollary 5.1.5 and Proposition 5.1.8]{BONFIGLIOLI}.

The following result is a consequence of \cite[Corollary 5.4.5]{BONFIGLIOLI} (see also \cite[Corollary 1.6]{fol2}), which provides a condition for the local integrability of a power of $\vert \cdot \vert _{\mathcal{G}}$.

\begin{lemma}\label{integrabilidadegauge}
Let $C_0$ the mean value of $\vert \cdot \vert_{\mathcal{G}}^{-Q}$, where $\vert \cdot \vert _{\mathcal{G}}$ is defined in \eqref{hom-norm-G}. Then if $\gamma>-Q$ and $0<a<b<\infty$,
$$ \int_{ a \leq  \vert y  \vert_{\mathcal{G}} \leq b} \vert y  \vert^{\gamma}_{\mathcal{G}}  d y=
\dfrac{C_0}{\gamma+Q}\left(b^{\gamma + Q}-a^{\gamma + Q}\right).
$$
\end{lemma}

As a consequence of Lemma \ref{integrabilidadegauge}, we obtain a sharper decay estimate than \eqref{est-semigrupo-unimod} for the $\|\cdot\|_{\infty}$ norm of the semigroup for initial data that behaves like a power of $\vert \cdot \vert_{\mathcal{G}}$ in the exterior of a ball.

\begin{lemma}\label{more-sharp-est} 
Let $\mathcal{G}$ be a homogeneous Carnot group and let
$u_0$ be a nonnegative measurable function satisfying
$$u_0(x) \geq \tilde{C}\,  \vert x \vert^{\gamma}_{\mathcal{G}},
\quad
\text{for almost every } x \in \mathcal G
\text{ such that }
\vert x \vert_{\mathcal G}\ge R, $$
where $-Q<\gamma<0$, and $\tilde{C}>0$, $R>0$ are given constants. Then there exists a constant $C>0$ such that
\begin{equation}\label{sem-esti-dados-radiais}
       \|S_{\mathcal{G}}(t) u_0\|_{\infty} \geq C\, t^{\frac{\gamma}{2}}
   \end{equation}
for all $t>R^2$.
\end{lemma}
\begin{proof}
Using the positivity of the heat kernel $H_t(y^{-1} \circ x)$, the Gaussian estimates \eqref{est-ker1} together with \eqref{pol-grow-hom}, the triangle inequality for $\vert \cdot   \vert_{\mathcal{G}}$ and Lemma \ref{integrabilidadegauge}, we obtain
\begin{align*}
\begin{aligned}
        &\left( S_{\mathcal{G}}(t) u_0 \right) (x) \\
        \noalign{\medskip}
        &\geq 
    \tilde{C} \int_{\vert y  \vert_{\mathcal{G}}\geq R} H_t(y^{-1} \circ x) \,  \big\vert  y  \big\vert^{\gamma}_{\mathcal{G}} d y\\
    \noalign{\medskip}
    &\geq C_1 t^{-\frac{Q}{2}}
    \int_{\vert y  \vert_{\mathcal{G}}\geq R}  \exp \left( -\frac{C_2 d_C( y^{-1} \circ x )^2}{t}\right)  \,  \big\vert  y  \big\vert^{\gamma}_{\mathcal{G}} d y\\
    \noalign{\medskip}
        &\geq C_1 t^{-\frac{Q}{2}}
    \int_{\vert y  \vert_{\mathcal{G}}\geq R}  \exp \left( -\frac{C_3 \vert y^{-1} \circ x \vert^2_{\mathcal{G}}}{t}\right)  \,  \big\vert  y  \big\vert^{\gamma}_{\mathcal{G}} d y\\
    \noalign{\medskip}
   & \geq C_1 t^{-\frac{Q}{2}} \exp \left( -\frac{  2C_4\vert x \vert^2_{\mathcal{G}} }{t}\right) \int_{ \sqrt{t} \leq  \vert y  \vert_{\mathcal{G}} \leq 2\sqrt{t}}  \exp \left( -\frac{  2C_5\vert y \vert^2_{\mathcal{G}}   }{t}\right)  \big\vert  y  \big\vert^{\gamma}_{\mathcal{G}}  d y\\
   \noalign{\medskip}
   & \geq C_1 \exp(-8C_5) \exp \left( -\frac{  2C_4\vert x \vert^2_{\mathcal{G}} }{t}\right) t^{-\frac{Q}{2}} \int_{ \sqrt{t} \leq  \vert  y   \vert_{\mathcal{G}} \leq 2\sqrt{t}} \big\vert  y  \big\vert^{\gamma}_{\mathcal{G}}  d y\\
   \noalign{\medskip}
   & \geq C_1 \exp(-8C_5) \exp(-2C_4) t^{-\frac{Q}{2}} \int_{ \sqrt{t} \leq  \vert y  \vert_{\mathcal{G}} \leq 2\sqrt{t}} \big\vert  y  \big\vert^{\gamma}_{\mathcal{G}}  d y\\
   \noalign{\medskip}
   & = C_1 \exp(-8C_5) \exp(-2C_4) \, C_0 \, (\gamma+Q)^{-1} \, \big[2^{\gamma+Q}-1\big] \, t^{-\frac{Q}{2}} t^{\frac{\gamma+Q}{2}}\\
   \noalign{\medskip}
   & = C t^{\frac{\gamma}{2}}\\
\end{aligned}
\end{align*}
whenever $  \vert x  \vert_{\mathcal{G}} \, \leq  \sqrt{t}$.
\end{proof}

We now present a lemma, recently established in \cite[Theorem 2.4]{Freire1}, which yields a refinement of estimate \eqref{decay-1}.

\begin{lemma}\label{Lemachave3} For all $t>0$ and $-Q <\alpha\leq0$, there exists a constant $C_0>0$ such that
$$
    \left\|   S_{\mathcal{G}}(t) \vert \cdot \vert^{\alpha}_{ \mathcal{G} }  \right\|_{\infty} \leq  C_0\,t^{\frac{\alpha}{2}} .
$$
\end{lemma}

\subsection{Case $p^- \leq 1+2(q+1)/Q < p^+$}\label{crit-case-section}

As a direct consequence of Lemma~\ref{more-sharp-est}, we obtain the following nonexistence result for global solutions.

\begin{theorem}\label{theo-carn} Let $h(t)=t^q$ with $q> -1$, and let $p$ be a bounded measurable function satisfying \eqref{ExponentVar} with $p^- >1$. Set $p^*= 1 + \frac{2(q+1)}{Q}$, and assume that either $p^+ \geq  p^* \geq  p^-$ or  $p^+  \geq  p^-\geq  p^* $. Let $f \in C([0,\infty))$ be such that
$$
f(s) \geq K s,$$
for some constant $K>0$ and all $s\geq 0$.
Then, for every nonnegative measurable function $u_0$ satisfying
\begin{equation}\label{eq-ref1-cacrit}
    u_0(x) \geq \tilde{C}\,  \vert x \vert^{\gamma}_{\mathcal{G}},
\quad
\text{for almost every } x \in \mathcal G
\text{ such that }
|x|_{\mathcal G}\ge R,
\end{equation}
with $\tilde{C}>0$, $R>0$ and
\begin{equation}\label{hip-pr-casocrit}
0< -\gamma < \frac{2(q+1)}{p^+-1}  \leq  Q ,
\end{equation}
the problem \eqref{eq1} (with $\mathcal{G}=\mathcal{M}$) does not admit a global nonnegative solution $u$ that is finite for almost every $(x,t)\in \mathcal{G}\times(0,\infty)$.
\end{theorem}
\begin{proof} Note that conditions (a) and (b) of Theorem \ref{thm1} (ii) are readily verified with $g(s) = Ks$, and it remains only to verify condition~\eqref{TT2} in order to apply the conclusion of Theorem~\ref{thm1} (ii).

If $\tau>\max \{1, R^2\}$ is sufficiently large, then estimate \eqref{sem-esti-dados-radiais} of Lemma \ref{more-sharp-est} holds and $C\tau^{\frac{\gamma}{2}}<1$. Therefore,
\begin{align*}
    \left(\int_{0}^{\tau} \sigma^q d \sigma \right)^{-1} &\int_{ \| S_{\mathcal{\mathcal{G}}}(\tau)u_0 \|_{\infty} }^{\infty} \frac{d \sigma}{\min \{\sigma^{p^-}, \sigma^{p^+}\}}\\
    \noalign{\medskip}
    &\leq  \left(\int_{0}^{\tau} \sigma^q d \sigma \right)^{-1} \int_{ C \tau^{\gamma/2} }^{\infty} \frac{d \sigma}{\min \{\sigma^{p^-}, \sigma^{p^+}\}}\\
    \noalign{\medskip}
    & = \frac{\tau^{-(q+1)}}{(q+1)^{-1}} \left(  \int_{ C \tau^{\gamma/2} }^{1} \frac{d \sigma}{\sigma^{p^+}} + \int_{ 1}^{\infty} \frac{d \sigma}{ \sigma^{p^-}}\right)\\
    \noalign{\medskip}
    & = \frac{\tau^{-(q+1)}}{(q+1)^{-1}} \left(  \frac{C^{-p^+ +1} \tau^{-\frac{\gamma}{2}(p^+ -1)}}{p^+ -1}  - \frac{1}{p^+ -1} + \frac{1}{p^- -1}\right).
\end{align*}
Since $p^{+}<1-2(q+1)/\gamma$ (by \eqref{hip-pr-casocrit}), we have that  
$$ \lim_{\tau \to \infty}  \frac{\tau^{-(q+1)}}{(q+1)^{-1}} \left(  \frac{C^{-p^+ +1} \tau^{-\frac{\gamma}{2}(p^+ -1)}}{p^+ -1}  - \frac{1}{p^+ -1} + \frac{1}{p^- -1}\right) =0.$$ 
Therefore, it is possible to choose $\tau > 1$ sufficiently large such that  
$$  \frac{\tau^{-(q+1)}}{(q+1)^{-1}} \left(  \frac{C^{-p^+ +1} \tau^{-\frac{\gamma}{2}(p^+ -1)}}{p^+ -1}  - \frac{1}{p^+ -1} + \frac{1}{p^- -1}\right) < 2^{-p^+} .$$  
Thus, estimate \eqref{cond-verf} is verified. In view of Theorem \ref{thm1} (ii), problem \eqref{eq1} with initial condition $u_0$ does not admit a global solution.
\end{proof}

Consider the class
\begin{equation*}
   \mathcal{I}_{a}= \Big\{ \phi \in L^{\infty}(\mathcal{G}); \phi \geq 0 \text{ and }  \liminf_{ \vert x\vert _{\mathcal{G}  } \to  \infty  } \big\vert  x  \big\vert_{\mathcal{G}}^a \, \phi(x) >0 \Big\}
\end{equation*}
where $a\geq 0$. In the Euclidean setting $\mathcal{G}=(\mathbb{R}^{N},+)$, with $L^{\infty}(\mathbb{R}^{N})$ replaced by $C_{b}(\mathbb{R}^{N})$, the space of bounded continuous functions on $\mathbb{R}^{N}$, the class $\mathcal{I}_{a}$ was introduced by Lee and Ni in \cite{LeeNi} to characterize the temporal behavior of the $L^{\infty}$-norm of solutions to the linear heat equation; see \cite[Lemma 2.12]{LeeNi}.

Note that if $\gamma<0$ satisfies \eqref{hip-pr-casocrit} and $u_0 \in \mathcal{I}_{-\gamma}$, then $u_0$ satisfies \eqref{eq-ref1-cacrit} for some sufficiently large $R>0$. 
Therefore, under the remaining assumptions of Theorem \ref{theo-carn}, the space $\mathcal{I}_{-\gamma}$ characterizes a class of initial data for which problem \eqref{eq1} does not admit a global solution. More precisely, we have the following.
\begin{corollary}\label{co-caso-crit1} 
Let $h(t)=t^q$ for $q>-1$ and $f(u)=u$. If either
$$
p^- \leq 1+\frac{2(q+1)}{Q} \leq p^+   \mbox{ or }  1+\frac{2(q+1)}{Q}  \leq p^-  \leq p^+,
$$
and $0 \neq u_0 \in \mathcal{I}_{-\gamma}$ with $\gamma$ satisfying \eqref{hip-pr-casocrit}, 
then problem~\eqref{eq1} with initial data $u_0$ does not admit a global solution.
\end{corollary}

\begin{remark}\label{VariableCrit} Here are some comments on the above Corollary. 
\begin{itemize}
\item[(i)] Note that the result in Corollary \ref{co-caso-crit1} is particularly interesting, as it applies to any variable exponent $p(x)$ satisfying \eqref{ExponentVar}.

\item[(ii)] In the Euclidean setting $\mathcal{G}=\mathbb{R}^N$ \cite[Examples 3.6 and 3.7]{Ferreira}, Ferreira et al. investigated the case 
$p^- < 1+\tfrac{2}{N} < p^+$ (that is, $q=0$). They proved that the associated parabolic problem admits 
functions $p(x)$ for which all solutions blow up in finite time, as well as functions $p(x)$ 
for which nontrivial global solutions exist.

More precisely, in the latter case, global solutions are constructed by means of a stationary supersolution associated with a suitably chosen discontinuous exponent. More precisely, one considers an exponent taking the value $p^-$ inside a ball and $p^+$ outside, and constructs a radial supersolution by matching an explicit decaying solution in the exterior region with a solution of an elliptic problem in the interior. We attempted to extend this construction to the subelliptic setting, since Theorem \ref{thm1} does not apply in the regime $p^- \leq p^* \leq p^+$ to guarantee the existence of global solutions.
However, two main difficulties arise. First, the comparison principle required to dominate the parabolic solution by a stationary supersolution is not known to hold in this sub-Riemannian framework. Second, the Euclidean construction relies crucially on the existence of an explicit radial exterior profile of the form $$u(x)=C \vert  x  \vert^{-\frac{2}{p^+-1}},$$ which solves the corresponding radial elliptic equation outside a ball. In contrast, on Carnot groups the sub-Laplacian of a radial function involves a nontrivial weight depending on the horizontal gradient of the homogeneous norm. In particular, in the Heisenberg group $\mathbb{H}^N$, if $\rho = |x|_{\mathbb{H}^N}$ denotes the standard homogeneous norm, one has
$$
\Delta_{\mathbb{H}^N} u(\rho)
=
\big|\nabla_{\mathbb{H}^N}\rho \big|^2
\left(
\frac{\partial^2 u}{\partial \rho^2}
+
\frac{Q-1}{\rho}\frac{\partial u}{\partial \rho}
\right),
$$
(see e.g. \cite[p. 299]{Cutri}). However, in general $|\nabla_{\mathbb{H}^N}\rho|^2$ is not constant on level sets of $\rho$, which prevents the reduction of the exterior problem to a one-dimensional ODE and obstructs the construction of explicit radial barriers of power type.

To the best of our knowledge, it is not known whether these difficulties can be overcome within the sub-Riemannian framework (see \cite{Cutri, Garofalo, JL1,JL2}), and they prevent a direct extension of the Euclidean supersolution construction. This leaves open the question of whether global solutions may exist in the intermediate regime
$$p^- \leq p^* \leq p^+$$
within the broader sub-Riemannian setting.

\end{itemize}
\end{remark}    

\subsection{Global Existence for a class of initial data}\label{critcase}

Let
\begin{equation*}
   \mathcal{I}^{a}= \left\{ \phi \in L^{\infty}(\mathcal{G}); \phi \geq 0 \text{ and }  \limsup_{ \vert x\vert _{\mathcal{G}  } \to  \infty  } \big\vert  x \big\vert_{\mathcal{G}}^a \, \phi(x) < \infty \right\}
\end{equation*}
where $a \geq 0$. This class is also closely related to a space introduced by Lee and Ni in \cite{LeeNi} in the Euclidean setting to study the long-time behavior of solutions to the heat equation.

As a direct consequence of Theorem \ref{thm1} and Lemma \ref{Lemachave3}, we obtain the following global existence result.
\begin{theorem}\label{theoGlobalLee}
Let $h(t) = t^q$ with $q > -1$, and let $p$ be a bounded measurable function satisfying \eqref{ExponentVar} with $p^- > 1$. 
Assume further that $p^+ \geq p^- > p^* = 1 + \tfrac{2(q+1)}{Q}$, and let 
$f \in C([0,\infty))$ be a nondecreasing function such that
$$f(s) \leq K s, \quad s \in (0,s_0),$$
for some constant $K > 0$. Suppose $0 \leq u_0 \in L^{\infty}(\mathcal{G})$ with $u_0 > 0$ on a set of positive measure, and assume that
\begin{equation}\label{eqLee1}
u_0(x) \leq \tilde{C}\, \vert  x  \vert^{\gamma}_{\mathcal{G}}
\quad \text{for almost every } x \in \mathcal{G} \text{ with } \vert x \vert_{\mathcal{G}} \geq R,
\end{equation}
where $\tilde{C} > 0$ and $R > 0$, and
\begin{equation}\label{eqLee2}
0 < \frac{2(q+1)}{p^- - 1} < -\gamma < Q.
\end{equation}
Then there exists $\Lambda > 0$ such that the problem \eqref{eq1} (with $\mathcal{G} = \mathcal{M}$) admits 
a global nonnegative solution
$$u \in L^{\infty}\big([0,\infty), L^{\infty}(\mathcal{G})\big)$$
with initial data $\lambda u_0$, for all $\lambda \in (0,\Lambda)$.
\end{theorem}

\begin{proof} Let $0\leq u_0 \in L^\infty(\mathcal{G})$ satisfy \eqref{eqLee1}, with $\gamma$ fulfilling \eqref{eqLee2}. Note that we may assume that $u_0$ satisfies \eqref{eqLee1} for all $R>0$. Indeed, let $M = \|u_0\|_{\infty}>0$. Then, for $\vert x \vert_{\mathcal{G}}\, <R$, $$u_0(x) \leq M \leq \big(M\vert x \vert_{\mathcal{G}}^{-\gamma} \big) \vert x \vert_{\mathcal{G}}^{\gamma} \, \leq (MR^{-\gamma}) \vert x \vert_{\mathcal{G}}^{\gamma}.$$
Therefore, setting $C_1 = \max\{ MR^{-\gamma}, \tilde{C} \}$, where  $\tilde{C}$ is the constant appearing in  \eqref{eqLee1},  we obtain $u_0(x) \leq C_1 \vert x \vert_{\mathcal{G}}^{\gamma}$ for almost every  $x \in \mathcal{G}$.

Since $f(s)\leq Ks$ for $s \in (0,s_0)$, continuity implies that $f(0)=0$. Thus, by the nondecreasing property of $f$, we can choose $\lambda>0$ sufficiently small such that $\|\lambda u_0\|_{\infty}<s_0$ and
$$f\big(\|S_{\mathcal{G}}(\sigma) \lambda u_0\|_{\infty}\big)\leq f\big(\|\lambda u_0\|_{\infty}\big) \leq 1,$$
for all $\sigma>0$. This and Lemma \ref{Lemachave3} implies
\begin{eqnarray*}
& & \int_{0}^{\infty} h(\sigma) \mathcal{F} \left(\|S_{\mathcal{G}}(\sigma)\lambda u_0\|_{\infty} \right) d\sigma \\
\noalign{\medskip}
&=&  \int_{0}^{\infty} h(\sigma) \, \frac{\max \left\{f(\|S_{\mathcal{G}}(\sigma) \lambda u_0\|_{\infty})^{p^-},\,f(\|S_{\mathcal{G}}(\sigma) \lambda u_0\|_{\infty})^{p^+} \right\}}{\|S_{\mathcal{G}}(\sigma)\lambda u_0\|_{\infty}} d\sigma  \\
\noalign{\medskip}
& = & \int_{0}^{\infty} \sigma^{q} \, \, \frac{f\big(\|S_{\mathcal{G}}(\sigma) \lambda u_0\|_{\infty}\big)^{p^-} }{\|S_{\mathcal{G}}(\sigma) \lambda u_0\|_{\infty}} d\sigma \\
\noalign{\medskip}
&\leq & K^{p^{-}} \lambda^{p^{-} -1}\int_{0}^{\infty} \sigma^{q} \, \, \|S_{\mathcal{G}}(\sigma)  u_0\|_{\infty}^{p^{-} - 1}  d\sigma \\
\noalign{\medskip}
&\leq & K^{p^{-}} \lambda^{p^{-} -1} \left( \frac{\| u_0\|_{\infty}^{p^{-} - 1}}{q+1} + \tilde{C} C_0\int_{1}^{\infty} \sigma^{q} \, \, \sigma^{\frac{\gamma}{2}(p^{-} - 1)}  d\sigma  \right).\\
 \end{eqnarray*} 
Condition \eqref{eqLee2} implies $\int_{1}^{\infty} \sigma^{q} \, \, \sigma^{\frac{\gamma}{2}(p^{-} - 1)}  d\sigma <\infty$, and thus the result follows from Theorem \ref{thm1}, by choosing $\lambda \in (0, \Lambda)$, where
$$ \Lambda = K^{^{-\frac{p^{-}}{p^- -1}}} \left( \frac{\| u_0\|_{\infty}^{p^{-} - 1}}{q+1} + \tilde{C} C_0\int_{1}^{\infty} \sigma^{q} \, \, \sigma^{\frac{\gamma}{2}(p^{-} - 1)}  d\sigma  \right)^{-\frac{1}{p^- -1}}.$$
\end{proof}
As a consequence of Theorem \ref{theoGlobalLee}, we obtain:
\begin{corollary}\label{co-caso-crit2} 
Let $h(t)=t^q$ with $q>-1$, and let $f \in C([0,\infty))$ be nondecreasing with 
$$
f(s) \leq K s \quad \text{for all } s \in (0,s_0),
$$
for some constant $K>0$. If 
$$
1+\frac{2(q+1)}{Q} < p^- \leq p^+,
$$
and $0 \neq u_0 \in \mathcal{I}^{-\gamma}$, where $\gamma$ satisfies \eqref{eqLee2}, then there exists $\Lambda>0$ such that the solution to problem~\eqref{eq1} with $\mathcal{G}=\mathcal{M}$ and initial data $\lambda u_0$ is global for all $\lambda \in (0,\Lambda)$.
\end{corollary}

\subsection{Second critical exponent in the sense of Lee and Ni \cite{LeeNi}} 
The following result is a direct consequence of Corollaries~\ref{co-caso-crit1} and~\ref{co-caso-crit2}.
It generalizes what is referred to in the Euclidean literature as the second critical Fujita exponent, 
a study initially developed by Lee and Ni in \cite{LeeNi} for the case $p = \text{constant}$ and $f(u)=u$. 
Since then, several contributions have extended this result to various Euclidean contexts; 
see, for example, \cite{WangJ1}, \cite{Tobakhanov}, and related works. 
Nevertheless, to our knowledge, results of this type are not yet available in the non-Euclidean setting.
\begin{corollary}\label{co-caso-crit3} 
Let $\mathcal{G}=\mathcal{M}$ be a homogeneous Carnot group, $h(t)=t^q$ with $q>-1$, and let $f \in C([0,\infty))$ be nondecreasing with 
$$
K_2 s \leq f(s) \leq K_1 s \quad \text{for all } s \in [0,\infty),
$$
for some constants $K_1, K_2>0$, and $1+\frac{2(q+1)}{Q} < p^- \leq p^+$.
\begin{itemize}
\item[(i)] If $0< -\gamma < \frac{2(q+1)}{p^+-1}  \leq  Q$ and $0\neq u_0 \in \mathcal{I}_{-\gamma}$, then problem~\eqref{eq1} subject to the initial data $u_0$, does not admit a global solution.
\item[(ii)] If $0 < \frac{2(q+1)}{p^- - 1} < -\gamma < Q$ and $0 \neq u_0 \in \mathcal{I}^{-\gamma}$, then there exists $\Lambda>0$ such that the solution to problem~\eqref{eq1} with initial data $\lambda u_0$ is global for all $\lambda \in (0,\Lambda)$.
\end{itemize}
\end{corollary}
\begin{remark} Here are some comments on the above Corollary
\begin{itemize}
\item[(i)] In our setting the values $\tfrac{2(q+1)}{p^+ - 1}$ and $\tfrac{2(q+1)}{p^- - 1}$ are the second critical Fujita exponents in the Lee and Ni sense.
\item[(ii)] Note that the above corollary reveals an interesting gap that will be addressed in future work. 
Specifically, it does not cover global or nonglobal results when 
$$-\gamma \in \Big( \tfrac{2(q+1)}{p^+ - 1}, \tfrac{2(q+1)}{p^- - 1} \Big).$$ 
This suggests the existence of new critical sets, analogous to those considered in the study of the second critical Fujita exponent for constant $p$, as introduced by Lee and Ni \cite{LeeNi} in the Euclidean setting.
\end{itemize}
\end{remark}

\section*{Acknowledgment}
R. Freire is partially supported by UESB under grant No. 027/2025. M. Loayza is supported by CNPq - 313382/2023-9 and CAPES/COFECUB 88887.879175/2023-00. 
\section*{Data Availability}

Data sharing is not applicable to this article as no datasets were generated or analyzed during the current study.

\section*{Statements and Declarations}

The authors declare that they have no competing interests.

\end{document}